\documentclass[twoside,11pt]{article}
\title{} \author{} \date{}
\usepackage{latexsym,amssymb,times}
\input amssym.def
\newtheorem{te}{Theorem}[section]
\newtheorem{prop}[te]{Proposition}

\newtheorem{fac}[te]{Fact}

\newtheorem{cla}[te]{Claim}

\newtheorem{ex}[te]{Example}

\def\dok{\noindent{\bf Proof. }}
\def\kdok{\hfill $\Box$ \par \vspace*{2mm} }
\def\a{\alpha}
\def\b{\beta}
\def\g{\gamma}

\def\f{\varphi}
\def\p{\psi}
\def\o{\omega}
\def\k{\kappa}
\def\r{\rho}
\def\s{\sigma}

\def\P{{\mathbb P}}
\def\Q{{\mathbb Q}}
\def\B{{\mathbb B}}
\def\N{{\mathbb N}}
\def\X{{\mathbb X}}
\def\Y{{\mathbb Y}}

\def\A{{\mathbb A}}

\def\BG{{\mathbb G}}

\def\BL{{\mathbb L}}
\def\CP{{\mathcal P}}

\def\CC{{\mathcal C}}
\def\CD{{\mathcal D}}

\def\CT{{\mathcal T}}

\def\CN{{\mathcal N}}
\def\F{{\mathcal F}}

\def\la{\langle}
\def\ra{\rangle}
\def\dom{\mathop{\mathrm{dom}}\nolimits}
\def\ran{\mathop{\mathrm{ran}}\nolimits}

\def\Aut{\mathop{\rm Aut}\nolimits}
\def\Cond{\mathop{\rm Cond}\nolimits}

\def\PC{\mathop{\mathrm{PC}}\nolimits}
\def\Bad{\mathop{\rm Bad}\nolimits}
\def\ar{\mathop{\rm ar}\nolimits}
\def\Var{\mathop{\rm Var}\nolimits}
\def\At{\mathop{\mathrm{At}}\nolimits}
\def\Form{\mathop{\rm Form}\nolimits}
\def\Fv{\mathop{\mathrm{Fv}}\nolimits}
\def\Sent{\mathop{\mathrm{Sent}}\nolimits}
\def\Th{\mathop{\rm Th}\nolimits}
\def\qr{\mathop{\mathrm{qr}}\nolimits}
\def\Mod{\mathop{\rm Mod}\nolimits}

\def\Lim{\mathop{\rm {Lim}}\nolimits}
\def\EF{\mathop{\mathrm{EF}}\nolimits}

\def\Col{\mathop{\rm Col}\nolimits}

\begin{document}
\thispagestyle{plain}
\begin{center}
           {\large \bf {\uppercase{Back and forth systems of condensations}}}
\end{center}
\begin{center}
{\bf Milo\v s S.\ Kurili\'c}\footnote{Department of Mathematics and Informatics, University of Novi Sad,
                                      Trg Dositeja Obradovi\'ca 4, 21000 Novi Sad, Serbia,
                                      e-mail: milos@dmi.uns.ac.rs}
\end{center}
\begin{abstract}
\noindent
If $L$ is a relational language, an $L$-structure $\X$ is condensable to an $L$-structure $\Y$, we write $\X \preccurlyeq _c \Y$,
iff there is a bijective homomorphism (condensation) from $\X$ onto $\Y$.
We characterize the preorder $\preccurlyeq _c$, the corresponding equivalence relation of bi-condensability, $\X \sim _c \Y$,
and the reversibility of $L$-structures in terms of back and forth systems and the corresponding games.
In a similar way we characterize the $\CP_{\infty \o}$-equivalence
(which is equivalent to the generic bi-condensability) and the $\CP$-elementary equivalence of $L$-structures,
obtaining analogues of Karp's theorem and the theorems of  Ehrenfeucht and Fra\"{\i}ss\'{e}.
In addition, we establish a hierarchy between the similarities of structures considered in the paper.
Applying these results we show that
homogeneous universal posets are not reversible
and that a countable $L$-structure $\X$ is weakly reversible
(that is, satisfies the Cantor-Schr\"{o}der-Bernstein property for condensations)
iff its $\CP_{\infty \o}\cup \CN_{\infty \o}$-theory is countably categorical.

{\sl 2010 MSC}:
03C07, 
03C75, 
03C50, 
03E40, 
06A06. 

{\sl Key words}:
condensation, bi-condensability, reversibility, back and forth, Karp's theorem, Ehrenfeucht--Fra\"{\i}ss\'{e} games, infinitary languages
\end{abstract}
\section{Introduction}\label{S1}
In this paper we continue the investigation of the {\it condensational preorder} $\preccurlyeq _c$ on the class $\Mod _L$ of structures of a relational language $L$,
defined by $\X \preccurlyeq _c \Y$ iff there exists a bijective homomorphism (condensation) from $\X$ onto $\Y$.
We also consider some naturally related relations and properties:
first, the equivalence relation of {\it bi-condensability}, defined by $\X \sim _c \Y $ iff
$\X \preccurlyeq _c \Y$ and $\Y \preccurlyeq _c \X$,
second, the structures $\X\in \Mod _L$ with the property that $\Y \sim _c \X$ implies $\Y \cong \X$, for all $\Y\in \Mod _L$,
(i.e., satisfying the Cantor-Schr\"{o}der-Bernstein property for condensations)
called {\it weakly reversible}, and, in particular, the {\it reversible} structures (that is, the structures $\X$ having the property that
each self-condensation of $\X$ is an automorphism).

At first sight, the relations between structures and the properties of structures mentioned above are more of set-theoretical
than of model-theoretical character; for example, reversibility is not preserved under
bi-definability and elementary equivalence \cite{KDef,KuMo1}. But, on the other hand,
reversibility is an invariant of some forms of bi-interpretability \cite{KRet},
and extreme elements of classes of structures definable by some $L_{\infty \o}$ sentences,
as well as the structures simply definable in linear orders are reversible \cite{KuMo2,KDef}).

It turns out that, in the investigation of condensability and related phenomena,
restricting our consideration to the class of sentences naturally corresponding to
condensations we obtain a possibility to use several basic concepts and methods
of model theory. So, in Section \ref{S4}, modifying (essentially Cantor's) theorem saying that back-and-forth equivalent countable $L$-structures are isomorphic,
we characterize (bi-)condensability of  structures of the same size $\k \geq \o$ in terms of back and forth systems of condensations and the corresponding games.
In addition we characterize reversible structures in this way and, as an application,
show that homogeneous-universal posets (and, in particular, the countable random poset) are non-reversible structures.

The main statements of Sections \ref{S5} and \ref{S6} are ``condensational analogues" of some well known results (concerning isomorphism).
Namely, first, it is evident that isomorphism of two $L$-structures implies that they satisfy the same $L_{\infty \o}$-sentences, which
implies their elementary equivalence.

Second, by the well known results including Karp's theorem,
the second property -- the $L_{\infty \o}$-equivalence of $L$-structures $\X$ and $\Y$,
their back and forth equivalence (partial isomorphism),
the existence of a winning strategy for player II in the corresponding Ehrenfeucht-Fra\"{\i}ss\'{e} game of length $\o$, $\EF _\o (\X ,\Y )$,
and the generic isomorphism of structures ($V[G]\models\X \cong \Y$, where $V[G]$ is some generic extension of the universe)
are equivalent conditions (see \cite{Karp}, \cite{Bar}, \cite{Nad}).

Third, by the classical results of Ehrenfeucht and Fra\"{\i}ss\'{e},
the third property -- the elementary equivalence,
the finitary isomorphism of $\X$ and $\Y$,
and the existence of winning strategies for player II in the games $\EF _n (\X ,\Y )$, for all $n\in \N$,
are equivalent conditions in the class of models of a finite language \cite{Ehr,Fra1,Fra2}.

So, roughly speaking, the results of Sections \ref{S5} and \ref{S6} show that,
replacing isomorphism by bi-condensability, $L_{\infty \o}$-equivalence by $\CP_{\infty \o}$-equivalence,
and elementary equivalence by $\CP$-equivalence, we obtain the analogues of all aforementioned classical theorems.
Of course, instead of back and forth systems of partial isomorphisms, back and forth systems of partial condensations come to the scene;
also, Ehrenfeucht-Fra\"{\i}ss\'{e} games are replaced by similar games with a different winning criterion.
As a by-product, in Section \ref{S5} we obtain the following  characterization: a countable structure $\X$ is weakly reversible iff the theory $\Th _{\CP_{\infty \o}\cup \CN_{\infty \o}}(\X )$ is $\o$-categorical.

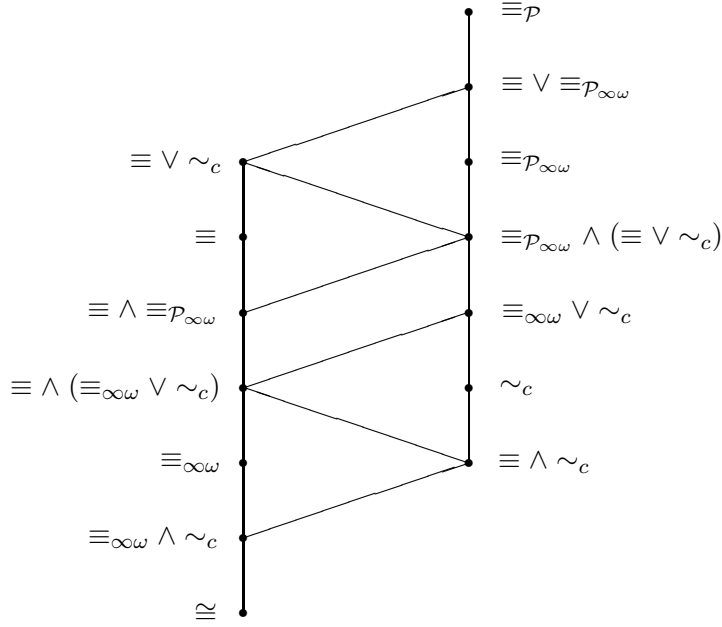
\begin{figure}
\unitlength 1mm 
\linethickness{0.5pt}
\ifx\plotpoint\undefined\newsavebox{\plotpoint}\fi 

\begin{picture}(110,90)(-14,0)


\put(35,5){\line(0,1){60}}
\put(65,25){\line(0,1){60}}
\put(35,15){\line(3,1){30}}
\put(65,25){\line(-3,1){30}}
\put(35,35){\line(3,1){30}}
\put(35,45){\line(3,1){30}}
\put(65,55){\line(-3,1){30}}
\put(35,65){\line(3,1){30}}

\put(35,5){\circle*{1}}
\put(35,15){\circle*{1}}
\put(35,25){\circle*{1}}
\put(35,35){\circle*{1}}
\put(35,45){\circle*{1}}
\put(35,55){\circle*{1}}
\put(35,65){\circle*{1}}
\put(65,25){\circle*{1}}
\put(65,35){\circle*{1}}
\put(65,45){\circle*{1}}
\put(65,55){\circle*{1}}
\put(65,65){\circle*{1}}
\put(65,75){\circle*{1}}
\put(65,85){\circle*{1}}

\put(30,5){\makebox(0,0)[cc]{$\cong  $}}
\put(23,15){\makebox(0,0)[cc]{$\equiv _{\infty \o}\land \sim _c  $}}
\put(28,25){\makebox(0,0)[cc]{$\equiv _{\infty \o}  $}}
\put(18,35){\makebox(0,0)[cc]{$\equiv \land \;(\equiv _{\infty \o}\lor \sim _c)  $}}
\put(23,45){\makebox(0,0)[cc]{$\equiv \land \equiv _{\CP _{\infty \o}} $}}
\put(30,55){\makebox(0,0)[cc]{$\equiv $}}
\put(26,65){\makebox(0,0)[cc]{$\equiv \lor \sim _c  $}}
\put(75,25){\makebox(0,0)[cc]{$\equiv \land \sim _c   $}}
\put(71,35){\makebox(0,0)[cc]{$\,\sim _c   $}}
\put(78,45){\makebox(0,0)[cc]{$\equiv _{\infty \o}\lor \sim _c  $}}
\put(84,55){\makebox(0,0)[cc]{$\equiv _{\CP _{\infty \o}}\land \;(\equiv \lor \sim _c) $}}
\put(74,65){\makebox(0,0)[cc]{$\equiv _{\CP _{\infty \o}}  $}}
\put(78,75){\makebox(0,0)[cc]{$\equiv \lor \equiv _{\CP _{\infty \o}} $}}
\put(72,85){\makebox(0,0)[cc]{$\equiv _{\CP }  $}}
\end{picture}

\vspace{-5mm}
\caption{Similarities of structures \label{F0003}}
\end{figure}

In Section \ref{S7} we compare the similarities of structures considered in this paper and show that the implications between them are as Figure \ref{F0003} describes.

We note that our restriction to relational structures is not essential. By \cite{KArb}
all results of this paper are in fact true for the structures of any language.

\section{Preliminaries}\label{S2}
\paragraph{Classes of formulas}
Let $L=\la R _i : i\in I\ra$ be a relational language, where $\ar (R_i)=n_i$, for $i\in I$,
let $\k$ be an infinite cardinal and $\Var =\{ v_\a : \a \in \k\}$ a set of variables.
By $\At _{L }$ we denote the corresponding set of  {\it atomic formulas}, that is,
$\At _{L } =  \{ v_\a =v_\b : \a ,\b  <\k\} \cup \{ R_i (v_{\a _1}, \dots , v_{\a_{n_i}}): i\in I \land  \la \a_1, \dots , \a_{n_i} \ra \in \k ^{n_i} \}$.

We recall that the class $\Form _{L_{\infty  \o }}$ of {\it $L_{\infty  \o }$-formulas} is the closure of the set $\At _{L }$
under negation, arbitrary conjunctions and disjunctions and finite quantification
(that is, $\neg \, \f$, $\forall v \,\f $ and $\exists v \,\f $ are in $\Form _{L_{\infty  \o }}$, whenever $\f \in \Form _{L_{\infty  \o }}$,
and $\bigwedge \F $ and $\bigvee \F  $ are in $\Form _{L_{\infty  \o }}$, for each set $\F  \!\subset\! \Form _{L_{\infty  \o }}$).

The set $\CP$ of {\it $R$-positive first order $L$-formulas} is the closure of the set
$$
\CP _0  =  \At _{L }  \cup \{ \neg \,v_\a = v_\b : \a ,\b<\k \}
$$
under finite conjunctions, disjunctions and quantification
(that is, $\f \land \p$, $\f \lor \p$, $\forall v \,\f $ and $\exists v \,\f $ are in $\CP $, whenever $\f ,\p \in \CP$;  negations are not allowed).

The set $\CN$ of {\it $R$-negative first order $L$-formulas} is the closure of the set
\begin{eqnarray*}
\CN _0 & := &  \{\neg\, R_i (v_{\a _1}, \dots , v_{\a_{n_i}}): i\in I \land  \la \a_1, \dots , \a_{n_i} \ra \in \k ^{n_i} \} \cup\\
       &    & \{ v_\a =v_\b : \a ,\b  <\k\} \cup \{ \neg \,v_\a = v_\b : \a ,\b<\k \}
\end{eqnarray*}
under finite conjunctions, disjunctions and quantification.
The class $\CP_{\infty  \o }$ (resp.\ $\CN _{\infty  \o }$) of {\it $R$-positive} (resp.\ {\it $R$-negative}) {\it $L_{\infty  \o }$-formulas}
is the closure of the set $\CP _0$ (resp.\ $\CN _0$)
under finite quantification and arbitrary conjunctions and disjunctions.

If $\F$ is a class of formulas, then $\Sent _\F$ will denote the class of all sentences from $\F$ and
$\F (v_0 , \dots ,v_{n-1})$ will be the set of formulas $\f \in \F $ such that
$\Fv (\f )\subset \{ v_0 , \dots ,v_{n-1}\}$.
If $\X ,\Y \in \Mod _L$,
by $\X \equiv _\F \Y$ (resp.\ $\X \lll _\F \Y$) we denote that
$\X \models \f $ iff $\Y \models \f $ (resp.\ $\X \models \f $ implies $\Y \models \f $), for each sentence $\f\in \Sent _\F$.

$L_{\infty  \o }$-formulas $\f$ and $\p$ are {\it logically equivalent}, in notation $\f \leftrightarrow \p$, iff
for each  $L$-structure $\X$ and any valuation
$\vec x \in {}^{\k } X$ we have: $\X \models \f [\vec x ]$ iff $\X \models \p [\vec x]$.
\paragraph{Back and forth systems of condensations}
If $\X $ and $\Y $ are $L$-structures,
a function $f$ will be called a {\it partial condensation from $\X$ to $\Y$}, we will write $f \in \PC (\X ,\Y )$, iff
$f$ is a bijection which maps $\dom f \subset X$ onto $\ran f \subset Y$ and
\begin{equation}\label{EQ0092}
\forall i\in I \;\; \forall \bar x \in (\dom f )^{n_i} \;\;(\bar x \in R_i^\X \Rightarrow f\bar x \in R_i^\Y ).
\end{equation}
\begin{fac}\label{T8198'}
If $n\in \N$, $ f =\{ \la x_k ,y_k \ra : k<n\}\subset X\times Y$, $\bar x =\la x_0 ,\dots , x_{n-1}\ra$ and $\bar y =\la y_0 ,\dots , y_{n-1}\ra$,
then $ f  \in \PC (\X ,\Y)$ iff
for each formula $\f \in \CP _0(v_0, \dots , v_{n-1})$ we have
\begin{equation}\label{EQ0096'}
\X \models \f [\bar x] \Rightarrow \Y \models \f [\bar y].
\end{equation}
\end{fac}
\dok
Clearly, $f$ is a function (resp.\ an injection) iff (\ref{EQ0096'}) is true for all formulas $v_{j_1} =v_{j_2}$ (resp.\ $\neg v_{j_1} =v_{j_2}$),
where $j_1,j_2 <n$, and $f$ satisfies (\ref{EQ0092}) iff (\ref{EQ0096'}) is true for all formulas
$R_i (v_{j_0}, \dots , v_{j_{n_i -1}} )$, where $i\in I$ and $j_k<n$, for $k<n_i$.
\kdok
A set $\Pi \subset \PC (\X , \Y )$ will be called a {\it back and forth system of condensations} (in the sequel, shortly: back and forth system, or b.f.s.)
iff $ \Pi\neq \emptyset$ and \\[-3mm]

(e1) $\forall f \in \Pi \;\, \forall x\in X  \; \exists g \in \Pi \;\; ( x\in \dom g \land f \subset g)$,

(e2) $\forall f \in \Pi \;\; \forall y\in Y \; \exists g \in \Pi \;\; ( y\in\, \ran g \,\;\land f \subset g )$.\\[-3mm]

\noindent
If such a b.f.s.\ exists we will write $\X  \preccurlyeq _{\rm c} ^{\mathrm{bfs}}\Y$ and
$\X  \sim _c ^{\mathrm{bfs}}\Y$ will denote that $\X  \preccurlyeq _{\rm c} ^{\mathrm{bfs}}\Y$ and $\Y  \preccurlyeq _{\rm c} ^{\mathrm{bfs}}\X$.
\paragraph{Games}
Let $\X$ and $\Y$ be $L$-structures and $\k$ a cardinal. {\it The game $G ^{\preccurlyeq _c} _{\k } (\X ,\Y )$}
is played in $\k$ steps by two players, I and II, in the following way: at the $\a$-th step player I chooses one of the two structures
and an element from it and then player II chooses an element from the other structure. More precisely, either
I chooses an $x_\a \in X$ and II chooses $y_\a \in Y$, or  I chooses an $y_\a \in Y$ and II chooses $x_\a \in X$.
So, each play gives a $\k$-sequence of pairs $\pi _\k =\la  \la x_\a ,y_\a \ra : \a <\k \ra \in {}^{\k}(X\times Y)$;
player II wins the play if for the set of pairs $\ran \pi _\k =\{ \la x_\a ,y_\a \ra : \a <\k \}\subset X\times Y$ we have
$\ran \pi _\k \in \PC (\X ,\Y)$. Otherwise, player I wins.

Roughly speaking, a {\it strategy} for a player determines its moves during the play on the basis of the previous moves of both players.
More formally, if $\Sigma $ is a strategy for player II, $\pi _\a =\la  \la x_\b ,y_\b \ra : \b <\a \ra$ is the sequence
produced in the first $\a$ moves and player I chooses $x_\a \in X$ at the $\a$-th step, then $y_\a :=\Sigma (\pi _\a ,x_\a) \in Y$ is the
response of II  suggested by $\Sigma$;  otherwise, if player I chooses $y_\a \in Y$, then $\Sigma$ suggests an $x_\a :=\Sigma (\pi _\a ,y_\a) \in X$.
A strategy $\Sigma $ is a {\it winning strategy} for a player iff  that player wins each play in which follows $\Sigma$.
We will write $\X  \preccurlyeq _c ^{G_{\k}}\Y$, if  player II has a winning strategy in the game $G_\k ^{\preccurlyeq _c}(\X ,\Y )$ and
$\X  \sim _c ^{G_{\k }}\Y$ will denote that $\X  \preccurlyeq _c ^{G_{\k }}\Y$ and $\Y  \preccurlyeq _c ^{G_{\k}}\X$.
\section{Condensability -- games of full length}\label{S4}
In this section we generalize the well known fact that back-and-forth equivalent countable $L$-structures are isomorphic.
We recall that a partial order $\P =\la P, \leq \ra $ is called {\it $\k$-closed} (where $\k$ is an infinite cardinal)
iff whenever $\g <\k$ is an ordinal and $\la p_\a : \a <\g \ra$ is a sequence in $P$ satisfying
\begin{equation}\label{EQ8300}
\forall \a ,\b \in \g \;\;(\a <\b \Rightarrow p_\b \leq p_\a )
\end{equation}
there is $p\in P$ such that $p\leq p_\a$, for all $\a <\g$.
If $\X ,\Y \in \Mod _L$, we will say that a b.f.s.\ $\Pi \subset \PC (\X ,\Y )$ is a {\it $\k$-closed b.f.s.}
iff the partial order $\la \Pi ,\supset\ra$ is $\k$-closed.
\begin{te}\label{T8195}
If $\X$ and $\Y$ are $L$-structures of size $\k\geq \o$, then we have

\noindent
(I) The following conditions are equivalent:

(a) $\X \preccurlyeq _c \Y$,

(b) There exists a $\k$-closed b.f.s.\ $\Pi\subset \PC (\X ,\Y )$,

(c) Player II has a winning strategy in the game $G ^{\preccurlyeq _c} _{\k }(\X ,\Y )$.

\noindent
(II) The following conditions are equivalent:

(a) $\X \sim _c \Y$,

(b) There are $\k$-closed b.f.s.\ $\Pi _{\X ,\Y}\subset \PC (\X ,\Y )$ and $\Pi _{\Y ,\X}\subset \PC (\Y ,\X )$,

(c) Player II has a winning strategy in the games $G ^{\preccurlyeq _c} _{\k }(\X ,\Y )$ and $G ^{\preccurlyeq _c} _{\k }(\Y ,\X )$.
\end{te}
A proof of the theorem is given after some preliminary work.
We recall that, if $\P =\la P, \leq \ra $ is a partial order, a set $D\subset P$ is called {\it dense} iff for each $p\in P$
there is $q\in D$ such that $q\leq p$. A set $\Phi\subset P$ is called a {\it filter} iff $\Phi \ni p\leq q$ implies $q\in \Phi$
and for each $p,q\in \Phi$  there is $r\in \Phi$ such that $r\leq p,q$.
\begin{fac}\label{T8192}
If $\,\P $ is a $\k$-closed partial order and $\CD$ a family of $\leq \k$ dense subsets of $\P$,
then there is a filter $\Phi$ in $\P$ intersecting all $D\in \CD$.
\end{fac}
\dok
Let $\CD =\{ D_\a :\a <\k \}$ be an enumeration.
By recursion we construct a sequence $\la p_\a :\a <\k \ra$ in $P$ such that for all $\a ,\b \in \k$ we have
(i) $p_\a \in D_\a$, and
(ii) $\a <\b \Rightarrow p_\b \leq p_\a $.
First we take $p_0\in D_0$. Suppose that $0<\a <\k$ and that $\la p_\b :\b <\a \ra$ is a sequence satisfying (i) and (ii).
Then, since $\P $ is $\k$-closed, there is $p\in P$ such that $p\leq p_\b$, for all $\b <\a$, and, since $D_\a$ is dense in $\P$,
there is $p_\a \in D_\a$, such that $p_\a \leq p$, which implies that $p_\a \leq p_\b$, for all $\b <\a$.
Thus the sequence $\la p_\b :\b \leq\a \ra$ satisfies (i) and (ii) and the recursion works.
It is evident that $\Phi :=\{ p\in P : \exists \a <\k \; p_\a \leq p\}$ is a filter  in $\P$.
\kdok
\begin{prop}\label{T8190}
If $\X$ and $\Y$ are $L$-structures of size $\k \geq \o $ and $\,\Pi\subset \PC (\X ,\Y )$ is a $\k$-closed b.f.s., then
each $f\in\Pi$ extends to a condensation $F\in \Cond (\X ,\Y )$.
\end{prop}
\dok
It is evident that the poset $\P :=\la \Pi _f , \supset\ra$, where $\Pi _f := \{ g \in \Pi : f \subset g \}$, is $\k$-closed.
Let $A:=X\setminus \dom f$ and $B:=Y\setminus \ran f$.

For $a \in A$, let $D_a :=\{ g \in \Pi _f : a \in \dom g \}$. If $ h  \in \Pi _f \setminus D_a$, then  $a \not\in \dom  h $
and by (e1) there is $g \in \Pi$ such that $a\in \dom g$ and $ h  \subset g $. Thus, since $f \subset  h $ we have $g\in \Pi _f$
and, hence, $g \in D_a$ and $g \supset  h $. So the sets $D_a$, $a\in A$, are dense in $\P$ and, similarly,
the sets $\Delta _b :=\{ g \in \Pi _f : b \in \ran g \}$, $b\in B$, are dense in $\P$ as well.
Since $|A|+|B|\leq \k$, by Fact \ref{T8192} there is a filter $\Phi $ in $\P$ intersecting all $D_a$'s
and $\Delta _b$'s. Clearly we have $f \subset F:= \bigcup \Phi\subset X \times Y $.

If $\la x,y'\ra , \la x,y''\ra \in F$ there are $g ' ,g '' \in \Phi$ such that $\la x,y'\ra \in g'$ and $\la x,y''\ra \in g ''$
and, since $\Phi$ is a filter, there is $g \in \Phi$ such that $g \supset g ' , g ''$. Thus $\la x,y'\ra , \la x,y''\ra \in g$
and, since $g$ is a function, $y'=y''$. So $F$ is a function and in the same way we show that it is an injection.

Let $i\in I$, $\bar x =\la x_0, \dots ,x_{n_i-1}\ra\in (\dom F )^{n_i}$ and $\bar x \in R_i^\X$.
Then, since $\dom F=\bigcup _{g \in \Phi}\dom g$, for each $j<n_i$ there is $g _j\in \Phi$ such that $x_j \in \dom g _j$.
Since $\Phi$ is a filter, there is $g \in \Phi$ such that $g \supset g _j$ and, hence, $\dom g \supset \dom g _j$, for all $j<n_i$.
Thus $\bar x\in (\dom g )^{n_i}$ and, since $g \in \PC (\X ,\Y )$ and $\bar x \in R_i^\X$, we have $F\bar x =g \bar x\in R_i^\Y$.
So  $F\in \PC (\X ,\Y )$.

If $a\in A$, then there is $g\in \Phi \cap D_a$ and, hence, $a\in \dom g \subset \dom F$; so  $A\subset \dom F$,
which, together with $\dom f\subset \dom F$, implies that $\dom F=X$.
Similarly we have  $\ran F=Y$ and, thus, $F\in \Cond (\X ,\Y )$.
\kdok
\noindent
{\bf Proof of Theorem \ref{T8195}}
We prove part (I), which evidently implies part (II).

(a) $\Rightarrow$ (b).
If $\X \preccurlyeq _c \Y$ and $F\in \Cond (\X ,\Y )$, then $\Pi := \{ F \} \subset \PC (\X , \Y )$.
Since $\dom F=X$ and $\ran F=Y$, the set $\Pi$ satisfies (e1) and (e2) trivially and
$\Pi$ is $\k$-closed because each sequence in $\Pi$ is a constant sequence.

(b) $\Rightarrow$ (a). If $\Pi\subset \PC (\X ,\Y )$ is a $\k$-closed b.f.s., then by Proposition \ref{T8190} there is $F\in \Cond (\X ,\Y )$.

(a) $\Rightarrow$ (c). Let $F\in \Cond (\X ,\Y )$.
Let $\Sigma$ be the following strategy for player II in the game $G ^{\preccurlyeq _c} _{\k }(\X ,\Y )$. At the $\a$-th step,
if player I chooses an $x_\a \in X$, then $\Sigma$ suggests $y_\a = F(x _\a)$;
if I chooses $y_\a \in Y$, then $\Sigma$ suggests $x_\a = F^{-1}(y _\a)$.
Now, if $\pi _\k=\la  \la x_\a ,y_\a \ra : \a <\k \ra $ is a play of the game in which player II follows $\Sigma$,
then $\ran \pi _\k =\{\la x_\a ,y_\a \ra : \a <\k \} \subset F$ and, hence, $\ran \pi _\k \in \PC (\X ,\Y )$; thus II wins the play.
So, $\Sigma $ is a winning strategy for player II.

(c) $\Rightarrow$ (a).
Let $\Sigma $ be a winning strategy for player II in the game $G ^{\preccurlyeq _c} _{\k }(\X ,\Y )$.
Clearly $\k = E\cup O$, where $E= \{ \g + 2n:\g \in \k\cap (\Lim \cup \{ 0\}) \land n\in \o \}$
and $O= \{ \g + 2n+1:\g \in \k\cap (\Lim \cup \{ 0\}) \land n\in \o \}$ are the sets of even and odd ordinals $<\k$ respectively, and we have
$|E|=|O|=\k$. Thus there is a bijection $b:\k \rightarrow X\cup Y$ such that $b[E]=X$ and $b[O]=Y$.
Let $\la  \la x_\a ,y_\a \ra : \a <\k \ra $ be the play of the game in which, at the step $\a$, player I chooses $b(\a )$
and player II follows $\Sigma$. Then $F=\{  \la x_\a ,y_\a \ra : \a <\k \}\in \PC (\X ,\Y)$.
If $x\in X$, then $x=b(\a )=x_\a$, for some $\a \in E$, and, hence, $x\in \dom F$; so $\dom F=X$ and, similarly,
$\ran F=Y$, which gives $F\in \Cond (\X ,\Y )$. Thus $\X \preccurlyeq _c \Y$ indeed.
\kdok
Proposition \ref{T8190} provides a useful characterization of reversible structures, which we prove in the sequel.
We note that the class of reversible structures contains, for example,
linear orders, Boolean lattices, well founded posets with finite levels \cite{Kuk},
tournaments, Henson graphs \cite{KuMo2}, and Henson digraphs \cite{KRet}.
Reversible equivalence relations are characterized in \cite{KuMo3}, while
reversible posets representable as disjoint unions of well orders and their inverses are characterized in \cite{KuMo4}.

If $\X$ is an $L$-structure, instead of $\PC (\X , \X)$ we write $\PC (\X )$.
A finite partial condensation $f \in \PC (\X )$ will be called {\it bad} iff $f \not\subset F$,
for all $F\in \Aut (\X )$.
\begin{te}\label{T8196}
For an $L$-structure $\X $ of size $\k\geq \o$ the following is equivalent:

(a) $\X$ is not a reversible structure,

(b) There exists a $\k$-closed b.f.s.\ $\Pi\subset \PC (\X )$ containing a bad condensation,

(c) There exists a b.f.s.\ $\Pi\subset \PC (\X )$ containing a bad condensation, if $\k =\o$.
\end{te}
\dok
(a) $\Rightarrow$ (b).
If $\X $ is not a reversible structure and $F\in \Cond (\X )\setminus \Aut (\X )$, then there are $i\in I$ and $\bar x =\la x_0 ,\dots ,x_{n_i-1}\ra\in X^{n_i}$
such that $\bar x \not\in R_i^\X$ and $F\bar x \in R_i^\X$.
Let $K:= \{ x_0 ,\dots ,x_{n_i-1} \}$ and $f  := F\upharpoonright K$.
Then, clearly, $f$ is a bad condensation, $f \in \Pi := \{ F\upharpoonright K : K\subset X \}\subset \PC (\X)$ and $\Pi$ is a $\k$-closed b.f.s.

(b) $\Rightarrow$ (a). If $\Pi\subset \PC (\X )$ is a $\k$-closed b.f.s.\ and $f\in \Pi$ a bad condensation,
then by Proposition \ref{T8190} there is $F\in \Cond (\X )$ extending $f$ and, hence, $F\not\in \Aut (\X )$. So, the structure $\X $ is not reversible.

The equivalence (b) $\Leftrightarrow$ (c) is true because each poset is $\o$-closed.
\kdok
\paragraph{Homogeneous universal posets}
Since the class of posets is a J\'{o}nsson class \cite{Jon}, for each regular beth number $\k$ there is a $\k$-homogeneous-universal poset $\P$.
As an example of application of Theorem \ref{T8196} we show that such posets are not reversible; in particular, taking $\k =\o$ we conclude that  the random poset
(i.e., the unique countable homogeneous universal poset, see \cite{Sch}) is non-reversible as well.
If $\mathbb P$ is a poset and $p,q \in P$, we will write $p\parallel q$ iff $p\not\leq q \land q\not\leq p$. For $A,B\subset P$,
$A<B$ denotes that $a<b$, for all $a\in A$ and $b\in B$; notation $A\parallel B$ is defined similarly.
\begin{te}\label{T8105}
Let $\P= \langle P, < \rangle$ be a strict partial order of size $\k \geq \o$ satisfying

{\rm (u1)} $\forall L,G \in [P ]^{<\k }\setminus \{ \emptyset \} \;\; (L<G \Rightarrow \exists x \in P \;\;L<x<G)$,

{\rm (u2)} $\forall K \in [P ]^{<\k }\setminus \{ \emptyset \} \;\; \exists x,y,z \in P \;\;(x<K \land y> K  \land z\parallel  K )$,

\noindent
and let $L,G,K \in [P ]^{<\k }\setminus \{ \emptyset \}$, where $L<G$. Then we have

(a) $| \{ x \in P : L<x<G \}|=\k$;

(b) $|\{ x \in P : x<K \}|=|\{ x \in P : x>K \}|=|\{ x \in P : x\parallel K \}|=\k$;

(c) If $\k$ is a regular cardinal, then $\P$ is not a reversible structure;

(d) If $\,\P$ is a $\k$-homogeneous-universal poset, it is not a reversible structure.
\end{te}
\dok
(a) If $L,G \in [P ]^{<\k }\setminus \{ \emptyset \}$, $L<G$ and $S:=\{ x \in P : L<x<G \}$, then by (u1) we have $S\neq \emptyset$ and, clearly,
$L<S$. So, assuming that $|S|<\k$, by (u1) there would be $y\in P$ such that $L<y<S$, which would imply that $y\not\in S$ and $L<y <G$.
But then $y\in S$ and we have a contradiction. So, $|S|=\k$.

(b) If $K \in [P ]^{<\k }\setminus \{ \emptyset \}$ then by (u1) we have  $T:=\{ x \in P : x<K \}\neq \emptyset$.
Assuming that $|T|<\k$, by (u1) there would be $x\in P$ such that $x<T$, which would imply that $x\not\in T$ and $x<K$.
But then $x\in T$ and we have a contradiction. So, $|T|=\k$ and, similarly, $|\{ x \in P : x>K \}|=|\{ x \in P : x\parallel K \}|=\k$.

(c) By (u2) 
there are $a_0,a_1,b_0,b_1\in P$, where $a_0 \parallel  a_1$ and $b_0 <b_1$.
Then $f _0:=\{ \la a_0,b_0\ra ,\la a_1,b_1\ra \}\in \Bad (\P )$, since $\k$ is a regular cardinal
the poset $\la \Pi \supset \ra$, where
$\Pi :=\{ f \in \PC  (\P ): f _0\subset f \land |f |<\k \} $, is $\k$-closed
and, by Theorem \ref{T8196}, it remains to be shown that $\Pi$ is a b.f.s.

(e1) If $f \in \Pi$ and $a\in P \setminus  \dom f$,
then the sets $L_a:=\{x\in \dom f : x<a\}$ and  $G_a:=\{y\in \dom f : y>a\}$ are of size $<\k$ and we have the following cases.

1. $L_a \neq \emptyset$ and $G_a \neq \emptyset$.
If $l\in f [L_a]$ and $g\in f [G_a]$,
then there are $x\in L_a$ and $y\in G_a$ such that $l=f (x)$ and $g=f(y)$
and, since $x<a<y$ and $f$ is a homomorphism, $f (x)<f (y)$, that is $l<g$.
Thus $f [L_a] < f [G_a] $
and, since $|\ran f|<\k$,
by (a) there is $b\in P \setminus \ran f$ such that $f [L_a] <b< f [G_a] $.
Now $g := f \cup \{ \la a,b \ra\}$ is an injection, $a\in \dom g$ and $f \subset g$.
If $x\in \dom f$ and $x<a$,
then $x\in L_a$
and, hence, $f (x)\in f [L_a]$ and $g(x)=f(x)<b =g (a)$.
Similarly, $a<y\in \dom f$ implies $g(a)< g(y)$.
So $g$ is a homomorphism and $g \in \Pi$.

2. $L_a = \emptyset$ and $G_a \neq \emptyset$.
Then $f [G_a]\in [P ]^{<\k }\setminus \{ \emptyset \}$
and, since $|\ran f|<\k$, by (b) there is $b\in P \setminus \ran f$ such that $b< f [G_a] $.
Now $g := f \cup \{ \la a,b \ra\}$ is an injection, $a\in \dom g$, $f \subset g$
and we show that $g$ is a homomorphism.
If $a<y\in \dom f$, then $y\in G_a$
and, hence, $f (y)\in f [G_a]$ and $g (a)=b <f (y)=g (y)$.

3. $L_a \neq \emptyset$ and $G_a = \emptyset$. This case is dual of case 2.

4. $L_a = \emptyset$ and $G_a = \emptyset$. Then $a\parallel  x$, for all $x\in \dom f$,
and choosing $b\in P\setminus \ran f$
we have $g := f \cup \{ \la a,b \ra\}\in \Pi$.

(e2) Let $f \in \Pi$ and $b\in D \setminus  \ran f$.
Since $|\dom f |<\k$ by (b) there is $a\in P\setminus \dom f$, such that  $a\parallel  x$, for all $x\in \dom f$.
Thus $g := f \cup \{ \la a,b \ra\} \in \Pi $ and (e2) is true indeed.

(d)
Let $\P$ be a $\k$-homogeneous-universal poset (of regular size $\k$).
Thus, each isomorphism between $<\k$-sized substructures of $\P$ extends to an automorphism of $\P$
and each poset of size $\leq \k$ embeds in $\P$. By (c), for a proof that $\P$ is not a reversible structure it is sufficient to
show that $\P$ satisfies (u1) and (u2).

Let $L,G \in [P ]^{<\k }\setminus \{ \emptyset \}$, where $L<G$, and let $\Y$ be a poset with domain $Y=Y_L \cup \{ a\} \cup Y_G$,
where $Y _L < \{ a\}< Y _G$, $\Y _L\cong \BL$
and  $\Y _G\cong \BG$. Then $|Y|<\k$ and, by the universality of $\P$ there is an embedding $e:\Y \hookrightarrow \P$.
Thus $L\cong e[Y_L]<e(a)<e[Y_B]\cong G$ and there are isomorphisms $f_L:e[Y_L]\rightarrow L$ and $f_G:e[Y_G]\rightarrow G$.
Clearly, $f:=f_L \cup f_G :e[Y_L]\cup e[Y_G]\rightarrow L\cup G$ is an isomorphism between $<\k$-sized substructures of $\P$
and, by the homogeneity of $\P$ there is an automorphism $F\in \Aut (\P )$ such that $f\subset F$, which implies that
$L <x:=F(e(a))<G$. So (u1) is true and (u2) has a similar proof.
\kdok
\section{Partial condensability -- games of countable length}\label{S5}
\begin{te}\label{T8194}
If $\X$ and $\Y$ are infinite  $L$-structures, then we have\\[-3mm]

\noindent
(I) The following conditions are equivalent:

(a) $\X \preccurlyeq _c \Y$, in some generic extension $V[G]$ of the universe,

(b) There is a b.f.s.\ $\Pi\subset \PC  (\X ,\Y )$,

(c) Player II has a winning strategy in the game $G ^{\preccurlyeq _c} _{\o }(\X ,\Y )$,

(d) $\X \lll _{\CP _{\infty \o}} \Y$,

(e) $\Y \lll _{\CN_{\infty \o}} \X$.

(f) $\X \preccurlyeq _c \Y$, if, in addition, the structures $\X$ and $\Y$ are countable.\\[-3mm]

\noindent
(II) The following conditions are equivalent:

(a) $\X \sim _c \Y$, in some generic extension $V_{\P }[G]$ of the universe,

(b) $\X  \sim _c ^{\mathrm{bfs}}\Y$, 

(c) $\X  \sim _c ^{G_\o}\Y$, 

(d) $\X \equiv _{\CP_{\infty \o}} \Y$,

(e) $\X \equiv _{\CN_{\infty \o}} \Y$,

(f) $\X \equiv _{\CP_{\infty \o} \cup \CN_{\infty \o}} \Y$,

(g) $\X \lll _{\CP _{\infty \o}\cup \CN_{\infty \o}} \Y$,

(h) $\X \sim _c \Y$, if, in addition, the structures $\X$ and $\Y$ are countable.
\end{te}
The proof is given after the following preliminaries.

To each $L_{\infty  \omega }$-formula $\varphi$ we adjoin a formula $\varphi  ^\neg $ in the following way.
First,
$(v_\alpha =v_\beta ) ^\neg  := \neg v_\alpha =v_\beta $ and
$(R_i (v_{\alpha _1}, \dots , v_{\alpha _{n_i}})) ^\neg  := \neg  R_i (v_{\alpha _1}, \dots , v_{\alpha _{n_i}}) $;
if $\varphi  ^\neg $ is defined for a formula $\varphi\in \Form _{L_{\infty  \omega }}$, then
$(\neg \varphi  ) ^\neg  := \varphi    $, $(\forall v_\alpha \; \varphi ) ^\neg  := \exists v_\alpha \; \varphi  ^\neg   $ and
$(\exists v_\alpha \; \varphi ) ^\neg  := \forall v_\alpha \; \varphi  ^\neg   $;
finally, if $\F  \subset \Form _{L_{\infty  \omega }}$ and $\varphi  ^\neg $ is defined for each formula $\varphi\in \F $, then
$(\bigwedge \F   ) ^\neg  :=\bigvee \F   ^\neg    $ and $(\bigvee \F   ) ^\neg  := \bigwedge \F   ^\neg  $,
where $\F   ^\neg  $ denotes the set $\{ \varphi  ^\neg  : \varphi \in \F  \}$.
The following statement is easily  provable by induction (see \cite{KuMo2}).
\begin{fac}\label{TA032}
Let $\varphi$ be an $L_{\infty  \omega }$-formula. Then

(a) $\varphi  ^\neg  \leftrightarrow \neg \varphi$;

(b) If $\varphi \in  \CP_{\infty  \omega }  $, then $\varphi  ^\neg \in  \CN_{\infty  \omega } $;

(c) If $\varphi \in \CN_{\infty  \omega } $, then $\varphi  ^\neg \in \CP_{\infty  \omega } $.
\end{fac}

\noindent
{\bf Proof of (I) of Theorem \ref{T8194}.}
(a) $\Rightarrow$ (b). Let $\P$ be a poset, $G$ a $\P$-generic filter over $V$, and let $F\in V_\P [G]$, where
$V_\P [G]\models F\in \Cond (\X ,\Y )$. Then $F=\tau _G$, for some $\P$-name $\tau$ and there is $p\in G$ such that
$p\Vdash ``\tau : \check{\X} \rightarrow  \check{\Y} \;\mbox{ is a condensation}"$.

We prove that
$\Pi :=\{ f \in \PC (\X ,\Y ):|f|<\o \land  \exists q\leq p \;\; q\Vdash \check{f}\subset \tau \}$
is a b.f.s. Let $f \in \Pi$, $a\in X\setminus \dom f$ and let $q\leq p$, where $q\Vdash \check{f } \subset \tau$.
Since $q\leq p$ and $p\Vdash \tau \in \Cond (\check{\X } ,\check{\Y } )$  we have $q\Vdash \tau \in \Cond (\check{\X } ,\check{\Y } )$
so, if  $H$ a $\P$-generic filter over $V$
and $q\in H$, in the extension $V_\P [H]$ we have: $\tau _H$ is a condensation from $\X$ to $\Y$ and $f \subset\tau _H$.
Now, for $b:=\tau _H (a)$ we have $f \subset g := f \cup \{ \la a,b \ra\}\subset \tau _H$ and, hence, there is $q_1\in H$
such that $q_1 \Vdash \check{g}\subset \tau$. Since $H$ is a filter and $q,q_1 \in H$, there is $q_2\in H$ such that $q_2 \leq q,q_1$.
Since $q_2 \leq q_1$ we have  $q_2 \Vdash \check{g}\subset \tau$ and $q_2 \leq q$ implies that $q_2\leq p $; thus $f \subset g \in \Pi$,
$a\in \dom g$ and (e1) is true. In a similar way we prove that $\Pi$ satisfies (e2).

(b) $\Rightarrow$ (c). Let $\Pi\subset \PC  (\X ,\Y )$ be a b.f.s. Since the closure of a b.f.s.\ under restrictions is a b.f.s.\ too,
w.l.o.g.\ we suppose that $\Pi$ is closed under restrictions. Let $<_X$ and $<_Y$ be well orders of the domains $X$ and $Y$ respectively and let us
pick $x^*\in X$ and $y^*\in Y$. Let $\Sigma$ be the strategy for player II in the game $G ^{\preccurlyeq _c} _{\o }(\X ,\Y )$ defined as follows.
Let $n\in \o$, let $\pi _n =\la \la x_k ,y_k\ra :k<n \ra$ be the sequence
produced in the first $n$ moves and $f_n =\{ \la x_k ,y_k\ra :k<n \}$. Now,

\noindent
-- if player I chooses $x_n \in X$ at the $n$-th step, then
$$
\Sigma (\pi _n ,x_n )=\left\{\begin{array}{cl}
                                 \min _{<_Y}\{ y\in Y : f _n \cup \{ \la x_n ,y \ra\}\in \Pi\} , &
                                                \mbox{if } f _n \cup \{ \la x_n ,y \ra\}\in \Pi, \\
                                                & \mbox{for some }y\in Y,\\[2mm]
                                 y^* , &  \mbox{otherwise;}
                            \end{array}
                     \right.
$$
-- if player I chooses $y_n \in Y$ at the $n$-th step, then
$$
\Sigma (\pi _n ,y_n )=\left\{\begin{array}{cl}
                                 \min _{<_X}\{ x\in X : f _n \cup \{ \la x ,y_n \ra\}\in \Pi\} , &
                                                \mbox{if } f _n \cup \{ \la x ,y_n \ra\}\in \Pi, \\
                                                & \mbox{for some }x\in X,\\[2mm]
                                 x^* , &  \mbox{otherwise.}
                            \end{array}
                     \right.
$$
Let $\pi _\o=\la \la x_k ,y_k\ra :k<\o \ra$ be a play in which player II follows $\Sigma$. By induction we show that
$f _n :=\{ \la x_k ,y_k\ra :k<n \}\in \Pi$, for all $n\in \o$.
First, $f _0=\emptyset \in \Pi$.
Suppose that $f _n\in \Pi$.
If, at the $n$-th step, player I have picked $x_n \in X$,
then by (e1) there is $g\in \Pi$ such that $f _n \subset g$ and $x_n\in \dom g$;
so $y:=g(x_n)\in Y$ and, since the set $\Pi$ is closed under restrictions,
$f _n \cup \{ \la x_n ,y \ra\}=g\upharpoonright (\dom f _n \cup \{ x_n \})\in \Pi$,
which implies that for $y_n=\Sigma (\pi _n ,x_n )$ we have $f _{n+1}=f _n \cup \{ \la x_n ,y _n\ra\}\in \Pi$.
If player I have picked $y_n \in Y$, then using (e2) we show that $f _{n+1}\in \Pi$ again.

So, $f _n \in \Pi \subset \PC (\X ,\Y )$, for all $n\in \o$, and, in addition, $f _0 \subset f_1 \subset \dots$ ,
which implies that $f =\bigcup _{n\in \o }f _n=\{ \la x_k ,y_k\ra :k<\o \}\in  \PC (\X ,\Y )$. Thus player II wins the play and $\Sigma$ is a winning strategy
for player II in the game $G ^{\preccurlyeq _c} _{\o }(\X ,\Y )$.

(c) $\Rightarrow$ (a). Let $\Sigma$ be a winning strategy
for player II in the game $G ^{\preccurlyeq _c} _{\o }(\X ,\Y )$
and let $\Pi$ be the set of the ranges $f_{\pi}=\{ \la x_k ,y_k\ra :k<n \}$ of all
finite partial plays  $\pi  =\la \la x_k ,y_k\ra :k<n \ra$, $n\in \o$, in which player II follows $\Sigma$.
We show that the sets $D_x :=\{ g \in \Pi : x\in \dom g \}$, $x\in X$, are dense in the poset $\P :=\la \Pi ,\supset \ra$.
So if $x\in X$ and $f _\pi  =\{ \la x_k ,y_k\ra :k<n \}\in \Pi$,
then, regarding a play in which $\pi  =\la \la x_k ,y_k\ra :k<n \ra$ is the sequence of the first $n$ moves,
$x$ can be considered as the choice $x_n \in X$ of player I at the $(n+1)$-st move.
Then player II takes $y_n=\Sigma (\pi ,x_n)$;
so we have $g := f_\pi \cup \{ \la x_n ,y_n \ra\}\in \Pi$
and $x\in \dom g$, thus $g\in D_x$ and $g \supset f_\pi$.
In a similar way we prove that the sets $\Delta _y :=\{ g \in \Pi : y\in \ran g \}$, $y\in Y$,
are dense in $\P$.
Now,  if $G$ is a $\P$-generic filter over $V$, then, since $\Pi \subset \PC (\X ,\Y )$,
in the generic extension $V_\P [G]$ we have $F:= \bigcup G \in \PC (\X ,\Y )$.
In addition, since $G\cap D_x\neq \emptyset$, for all $x\in X$, we have $\dom F =X$,
and, similarly, $\ran F =Y$; thus $F\in \Cond (\X ,\Y )$.

(b) $\Rightarrow$ (d). Let $\Pi \subset \PC (\X , \Y )$ be a b.f.s.
By induction on the construction of $\CP _{\infty \o}$-formulas we show that for each formula $\p (v_0, \dots ,v_{n-1})\in \CP _{\infty \o}$ we have
\begin{equation}\label{EQ8321}
\forall f\in \Pi \;\; \forall \bar x \in (\dom f)^n \;\;\Big( \X \models \p [\bar x]\Rightarrow\Y \models \p [f \bar x]\Big).
\end{equation}
Let $\p (\bar v)\in \CP _0$. If $f\in \Pi$, $\bar x \in (\dom f)^n$ and $\X \models \p [\bar x]$,
then $f$ is a condensation from
a substructure $\A$ of $\X$ onto
a substructure $\B$ of $\Y$ and $\bar x \in A^n$.
Since $\X \models \p [\bar x]$ and $\p \in \Sigma _0$ we have $\A \models \p [\bar x]$,
which, since $f$ is a condensation and $\p \in \CP _{\infty \o}$ implies $\B \models \p [f\bar x]$,
so, since $\p$ is a $\Sigma _0$-formula, $\Y \models \p [f\bar x]$ and (\ref{EQ8321}) is true.

Let $\p '(\bar v):= \bigwedge \F  $ and suppose that (\ref{EQ8321}) is true for all $\p (\bar v)\in \F $.
If $f \in \Pi $, $\bar x \in (\dom f)^n$ and $\X \models \p ' [\bar x]$,
then for each $\p \in \F $  we have $\X \models \p [\bar x]$ and, by (\ref{EQ8321}), $\Y \models \p [f \bar x]$,
which implies $\Y \models \p '[\bar y]$ and we are done.

Let $\p '(\bar v):= \bigvee \F  $ and suppose that (\ref{EQ8321}) is true for all $\p (\bar v)\in \F $.
If $f \in \Pi $, $\bar x \in (\dom f)^n$ and $\X \models \p ' [\bar x]$,
then for some $\p _0 \in \F $ we have $\X \models \p _0[\bar x]$ and, by (\ref{EQ8321}), $\Y \models \p _0 [f \bar x]$,
which implies that $\Y \models \p '[\bar y]$.

Let $\p '(\bar v):= \exists v_n \;\p (\bar v, v_n) $ and let (\ref{EQ8321}) be true for $\p (\bar v, v_n)$, that is
\begin{equation}\label{EQ8322}
\forall g\in \Pi \; \forall \bar x \in (\dom g)^n \;\forall x \in \dom g \;
\Big( \X \models \p [\bar x, x]\Rightarrow\Y \models \p [g\bar x ,g x]\Big).
\end{equation}
If $f \in \Pi $, $\bar x \in (\dom f)^n$ and $\X \models \p ' [\bar x]$, then there is $x \in X$
such that $\X \models \p [\bar x, x]$, by (e1) there is $g \in \Pi$ such that $x\in \dom g$ and $f \subset g$ and, by
(\ref{EQ8322}) and since $g\bar x =f\bar x$ we have $\Y \models \p [f\bar x ,gx]$, which gives  $\Y \models \p ' [f\bar x]$.

Let $\p '(\bar v):= \forall v_n \;\p (\bar v, v_n) $, and let (\ref{EQ8322}) be true for $\p (\bar v, v_n)$.
If $f \in \Pi $, $\bar x \in (\dom f)^n$ and $\X \models \p ' [\bar x]$,
then for each $x \in X$ we have $\X \models \p [\bar x, x]$.
For $y\in Y$  by (e2) there is $g \in \Pi$ such that $y\in \ran g$ and $f \subset g$;
thus $y=g(x)$, for some $x\in \dom g$. By our assumption we have $\X \models \p [\bar x, x]$
and by (\ref{EQ8322}), since $g\bar x =f\bar x$ we have $\Y \models \p [f\bar x ,y]$.
This gives  $\Y \models \p ' [f\bar x]$.

Now, by (\ref{EQ8321}), if $\p\in \Sent _{\CP _{\infty \o}}$ and $\X \models \p$, then $\Y \models \p$. Thus  $\X \lll _{\CP _{\infty \o}} \Y$.

(d) $\Rightarrow$ (b). Let $\X \lll _{\CP_{\infty \o}} \Y$. For $n<\o$, let $P_n$ be the set of pairs $\la \bar x ,\bar y \ra \in X^n \times Y^n$
such that
\begin{equation}\label{EQ8323}
\forall \p (v_0,\dots , v_{n-1})\in \CP _{\infty \o}\;\; \Big( \X \models \p [\bar x] \Rightarrow\Y \models \p [\bar y]\Big).
\end{equation}
Then, since $f_{\bar x , \bar y}:=\{ \la x_k , y_k\ra : k<n\}$ preserves all formulas of the form
$v_k =v_l$, $\neg\;v_k =v_l$ and $R_i(v_{k_0},\dots ,v_{k_{n_i -1}})$,
we have $\Pi :=\{ f_{\bar x ,\bar y}:\la \bar x ,\bar y \ra\in \bigcup _{n<\o}P_n \}\subset  \PC (\X ,\Y )$.
Since $\X \lll _{\CP_{\infty \o}} \Y$ we have $f_{\emptyset, \emptyset}=\emptyset \in \Pi$.

(e1) Let $f_{\bar x ,\bar y}\in \Pi$ and $a\in X\setminus \{ x_0 ,\dots ,x_{n-1}\}$. We show that there is $b\in Y $ such that
\begin{equation}\label{EQ8326}
 \forall \p (v_0,\dots , v_{n-1}, v_n)\in \CP_{\infty \o} \;\; \Big( \X \models \p [\bar x,a] \Rightarrow\Y \models \p [\bar y,b]\Big).
\end{equation}
On the contrary, there would be formulas $\p _b (v_0,\dots , v_{n-1}, v_n)$, $b\in Y$, such that
\begin{equation}\label{EQ8324}
\forall b\in Y \;\; \X \models \p _b[\bar x,a]
\end{equation}
\begin{equation}\label{EQ8325}
\forall b\in Y \;\; \Y \models \neg \p _b [\bar y,b]
\end{equation}
Now $\p (v_0,\dots , v_{n-1}):= \exists v_n \bigwedge _{b\in Y} \p _b (v_0,\dots , v_{n-1}, v_n)\in \CP_{\infty \o}$
and $\X \models \p [\bar x]$ iff there is $a\in X$ such that $\X \models \p _b[\bar x,a]$, for all $b\in Y$, which is true by (\ref{EQ8324}).
So, by (\ref{EQ8323}) we have $\Y \models \p [\bar y]$ and, hence, there is $b^*\in Y$ such that
$\Y \models \p _b[\bar y,b ^*]$, for all $b\in Y$. In particular, $\Y \models \p _{b^*}[\bar y,b^*]$,
which is false by (\ref{EQ8325}).

By (\ref{EQ8326}) we have $\la \bar xa, \bar yb\ra \in P_{n+1}$ so $f_{\bar x ,\bar y}\subset f_{\bar xa, \bar yb}\in \Pi$ and $a\in \dom f_{\bar xa, \bar yb} $.

(e2) Let $f_{\bar x ,\bar y}\in \Pi$ and $b\in Y\setminus \{ y_0 ,\dots ,y_{n-1}\}$. We show that there is $a\in X $ such that (\ref{EQ8326}) holds.
On the contrary, there would be formulas $\p _a (v_0,\dots , v_{n-1}, v_n)$, $a\in X$, such that
\begin{equation}\label{EQ8327}
\forall a\in X \;\; \X \models \p _a[\bar x,a]
\end{equation}
\begin{equation}\label{EQ8328}
\forall a\in X \;\; \Y \models \neg \p _a [\bar y,b]
\end{equation}
Now $\p (v_0,\dots , v_{n-1}):= \forall v_n \bigvee _{a\in X} \p _a (v_0,\dots , v_{n-1}, v_n)\in \CP_{\infty \o}$
and $\X \models \p [\bar x]$ iff  for each $x\in X$ there is $a\in X$ such that $\X \models \p _a[\bar x,x]$,
which is,  by (\ref{EQ8327}), true for $a=x$.
So $\X \models \p [\bar x]$ and by (\ref{EQ8323}) we have $\Y \models \p [\bar y]$.
Thus, for each $b^*\in Y$ and, in particular, for $b$, there is $a\in X$ such that $\Y \models \p _a[\bar y,b ]$,
which is false by (\ref{EQ8328}).

So, there is $a\in X $ such that (\ref{EQ8326}) holds
and we have $\la \bar xa, \bar yb\ra \in P_{n+1}$
so $f_{\bar x ,\bar y}\subset f_{\bar xa, \bar yb}\in \Pi$ and $b\in \ran f_{\bar xa, \bar yb} $.

(d) $\Leftrightarrow$ (e). Let $\X \lll _{\CP_{\infty \o}} \Y$, $\f \in \Sent _{\CN_{\infty \o}}$ and $\Y \models \f$. Assuming that $\X \models \neg \f$,
by Fact\ref{TA032}(a) and (c) we would have $\X \models  \f ^\neg$ and $\f ^\neg \in \CP_{\infty \o}$; so, since $\X \lll _{\CP_{\infty \o}} \Y$,
$\Y \models  \f ^\neg$ that is $\Y \models  \neg\f $, which gives a contradiction. So $\X \models  \f$ and, thus, $\Y \lll _{\CN_{\infty \o}} \X$.
The converse has a symmetric proof.

The equivalence (b) $\Leftrightarrow$ (f) for countable structures
follows from  part (I) of Theorem \ref{T8195} (the equivalence (a) $\Leftrightarrow$ (b)), because each poset is $\o$-closed.
\kdok
\noindent
{\bf Proof of (II) of Theorem \ref{T8194}.}
The equivalence of conditions (b), (c), (d) and (e) follows from (I).
If (d) holds, we have (e) as well and, hence (f) is true. The implication  (f) $\Rightarrow$ (g) is trivial.
If (g) is true, then $\X \lll _{\CP _{\infty \o}} \Y$ and $\X \lll _{\CN_{\infty \o}} \Y$, which by (I) implies that $\Y \lll _{\CP _{\infty \o}} \X$.
Thus $\X \equiv _{\CP _{\infty \o}}\Y$ and, thus (g) implies (d).

If (a) holds, then
$V_{\P }[G]\models \X \preccurlyeq _c \Y$ and $V_{\P }[G]\models \Y \preccurlyeq _c \X$; so by (I)
there are b.f.s.\ $\Pi _{\X ,\Y}\subset \PC (\X ,\Y )$ and $\Pi _{\Y ,\X}\subset \PC (\Y ,\X )$ and (b) is true.

If (b) is true,  $\P _0: =\la \Pi _{\X ,\Y} ,\supset \ra$, $\P _1: =\la \Pi _{\Y ,\X} ,\supset \ra$, $\P =\P_1 \times \P _2$ and $G$ is a
$\P$-generic filter over $V$,
then (see \cite{Kun}, p.\ 253) the generic extension $V_\P [G]$ contains
a  $\P_0$-generic filter over $V$, $G_0$,
and a  $\P_1$-generic filter over $V$, $G_1$, and, hence, $V_{\P _0} [G_0] ,V_{\P _1} [G_1] \subset V_\P [G]$.
Thus (see the proof of the implication (b) $\Rightarrow $ (a) of part (I)) the extension $V_\P [G]$
contains condensations $F: \X \rightarrow \Y$ and $G: \Y \rightarrow \X$, which means that $V_\P [G]\models \X \sim _c \Y$ and (a) is true as well.
\kdok
\begin{te}\label{T8102}
For a countable $\X \in \Mod _L$ the following conditions are equivalent:

(a) $\X$ is weakly reversible,

(b) $\X \equiv _{\CP_{\infty \o}}\! \Y $ implies  that $\Y\cong \X$, for each countable $\Y \in \Mod _L$,

(c) $\Th _{\CP _{\infty \o}\cup \CN_{\infty \o}} (\X ):=\{ \f \in \Sent _{\CP _{\infty \o}\cup \CN_{\infty \o}} : \X \models \f \}$ is $\o$-categorical.
\end{te}
\dok
(a) $\Leftrightarrow$ (b) follows from Theorem \ref{T8194} (part II, (d) $\Leftrightarrow$ (h)).
$\Th _{\CP _{\infty \o}\cup \CN_{\infty \o}} (\X )$
is an $\o$-categorical theory iff for each countable $\Y \in \Mod _L$, $\X \lll _{\CP_{\infty \o} \cup \CN_{\infty \o}}\! \Y $
that is, by Theorem \ref{T8194}, $\X \equiv _{\CP_{\infty \o}}\! \Y $, implies $\Y\cong \X$. So, (b) $\Leftrightarrow$ (c) is true.
\kdok
In addition, by Theorem \ref{T8194}, condition $\X \equiv _{\CP_{\infty \o}}\! \Y $ in (b)
can be replaced by $\X \equiv _{\CN_{\infty \o}}\! \Y $, $\X \equiv _{\CP_{\infty \o}\cup \CN_{\infty \o}}\! \Y $, or
$\X \lll _{\CP_{\infty \o} \cup \CN_{\infty \o}}\! \Y $.
\begin{ex}\label{EX0002}\rm
The theory $\Th _{\CP _{\infty \o}\cup \CN_{\infty \o}} (\X )$ is $\o$-categorical, but the first
order theory of $\X$, $\Th (\X )$, is not.
Let $\X$ be a countable structure with one equivalence relation, such that there are no infinite equivalence classes and
for each $n\in \N$ there is exactly one equivalence class of size $n$. It is known that $\Th (\X )$ is not $\o$-categorical.
But $\X$ is a reversible and, hence, a weakly reversible structure (see \cite{KuMo3} for a characterization of reversible equivalence relations).
So, by Theorem \ref{T8102}, the theory $\Th _{\CP_{\infty \o}\cup \CN_{\infty \o} } (\X )$ is $\o$-categorical.
\end{ex}
\section{Finitary condensability -- finite games}\label{S6}
Here we consider the relations $\lll _{\CP}$ and $\equiv _{\CP}$, when $L$ is a finite relational language.
We note that, if $L$ is an infinite language and $\X ,\Y \in \Mod _L$, then $\X \lll _{\CP}\Y$ iff
$\X |L'\lll _{\CP_{L'}}\Y|L'$, for each finite $L' \subset L$, and the same holds for the relation $\equiv _{\CP}$.

If $\X$ and $\Y$ are $L$-structures we will say that $\X$ is {\it finitely condensable} to $\Y$
iff  there is a sequence $\la \Pi _r : r<\o\ra$, where, for each $r<\o$,
$\emptyset \neq \Pi_r \subset \PC (\X ,\Y)$ and  \\[-2mm]

(f1) $\forall  f \in \Pi_{r+1} \;\; \forall x\in X \;\; \exists  g \in \Pi_r \;\; (x\in \dom  g  \land  f \subset  g )$,

(f2) $\forall  f \in \Pi_{r+1} \;\; \forall y\in Y \;\; \exists  g \in \Pi_r \;\; (y\in \,\ran  g  \;\land  f \subset  g )$.\\[-2mm]

\noindent
Then we will write $\X \preccurlyeq _c^{\mathrm{fin}} \Y$. If, in addition, $\Y \preccurlyeq _c^{\mathrm{fin}} \X$, we will write $\X \sim _c^{\mathrm{fin}} \Y$.
\begin{te}\label{T0099}
If $L$ is a finite language and $\X$ and $\Y$ are $L$-structures, then we have\\[-3mm]

\noindent
(I) The following conditions are equivalent:

(a) There is a sequence $\la \Pi _r : r<\o\ra$ satisfying (f1) and (f2), 

(b) For each $n\in \o$, player II has a winning strategy for $G_n^{\preccurlyeq _c}(\X ,\Y )$,

(c) $ \X \lll _\CP \Y $,

(d) $ \Y \lll _\CN \X $.\\[-3mm]

\noindent
(II) The following conditions are equivalent:

(a) $\X \sim _c^{\mathrm{fin}} \Y $,

(b) $\X  \sim _c ^{G_{n}}\Y$, for all $n\in \o$,

(c) $ \X \equiv _\CP \Y $,

(d) $ \X \equiv _\CN \Y $,

(e) $ \X \equiv _{\CP \cup \CN} \Y $,

(f) $ \X \lll _{\CP \cup \CN} \Y $.
\end{te}
\dok
We prove part I, which gives the equivalence of (a) -- (d) of part II. The proof of the rest is similar to the proof of the corresponding
part of Theorem \ref{T8194}.

(a) $\Rightarrow$ (c). Let $\X \preccurlyeq _c^{\mathrm{fin}} \Y$ and let $\la \Pi _r : r<\o\ra$ be a sequence which witnesses it.
By induction on the complexity of formulas we prove that
for each $\f (v_0, \dots , v_{n-1})\in \CP$  we have
\begin{equation}\label{EQ0093}
\forall r\geq \qr (\f )\;\; \forall  f \in \Pi_r \;\;\forall \bar x\in (\dom  f )^n\;\; \Big( \X \models \f [\bar x] \Rightarrow \Y \models \f [ f \bar x]\Big).
\end{equation}
If $\f \in \CP _0 $, then $\qr {\f }=0$.
So, if $r<\o$, $ f \in \Pi_r$, $\bar x\in (\dom  f )^n $ and $\X \models \f [\bar x]$, we show that $\Y \models \f [ f \bar x]$.
Since $ f $ is an injection, this is true for all formulas of the form $v_i =v_j$ and $v_i\neq v_j$.
If $\f := R_i (v_{j_0}, \dots , v_{j_{n_i -1}} )$, where $j_0, \dots , j_{n_i -1}< n_i$,
then $\la x_{j_0}, \dots , x_{j_{n_i -1}}\ra\in R_i ^\X$
and, by (\ref{EQ0092}),
$\la  f (x_{j_0}), \dots ,  f (x_{j_{n_i -1}})\ra\in R_i ^\Y$
which, since $ f \bar x=\la  f (x_0), \dots , f (x_{n_i -1})\ra$, gives $\Y \models \f [ f \bar x]$.

Suppose that (\ref{EQ0093}) holds for all formulas from $\CP _k$ and that $\f (v_0, \dots , v_{n-1})\in \CP _{k+1}$.
Let $r\geq \qr (\f )$, $ f \in \Pi_r $ and $\bar x\in (\dom  f )^n$.

If $\f = \p _1 \lor \p _2 $ or $\f = \p _1 \land \p _2 $,
then $\p _1 ,\p _2 \in \CP _k$, $\qr (\p _1),\qr (\p _2 )\leq \qr (\f )$ and, since (\ref{EQ0093}) holds for $\p_1$ and $\p_2$, we have
$$
\X \models \p _1 [\bar x] \Rightarrow \Y \models \p _1 [ f \bar x] \;\;\mbox{ and }\;\;
\X \models \p _2 [\bar x] \Rightarrow \Y \models \p _2 [ f \bar x] .
$$
Now, if $\f = \p _1 \lor \p _2 $ and $\X \models \f [\bar x]$,
then $\X \models \p _j [\bar x]$, for some $j\in 2$, which implies $\Y \models \p _j [ f \bar x]$;
and, hence, $\Y \models \f [ f \bar x]$.
If $\f = \p _1 \land \p _2 $ and $\X \models \f [\bar x]$,
then $\X \models \p _0 [\bar x]$ and $\X \models \p _1 [\bar x]$,
which implies $\Y \models \p _0 [ f \bar x]$ and $\Y \models \p _1 [ f \bar x]$
and, hence, $\Y \models \f [ f \bar x]$.

If $\f =\exists v_n \p (v_0, \dots , v_{n-1},v_n)$ and $\X \models \f [\bar x]$,
then there is $x\in X$ such that $\X \models \p [\bar x ,x]$,
by (f1) there is $ g \in \Pi_{r-1}$ such that $x\in \dom  g $ and $ f \subset  g $
and, since $\qr (\p)=\qr (\f )-1\leq r-1$,
by the induction hypothesis for $\p$, $r-1$, $ g $ and $\bar x ^\smallfrown x$,
from $\X \models \p [\bar x ,x]$ it follows that $\Y \models \p [ f \bar x , g (x)]$.
So there is $y= g (x)\in Y$ such that $\Y \models \p [ f \bar x ,y]$, which means that $\Y \models \f [ f \bar x]$.

If $\f =\forall v_n \p (v_0, \dots , v_{n-1},v_n)$ and $\X \models \f [\bar x]$, then
\begin{equation}\label{EQ0097}
\forall x\in X \;\;\X \models \p [\bar x ,x].
\end{equation}
If $y\in Y$,
then by (f2) there is $ g \in \Pi_{r-1}$ such that $y\in \ran  g $ and $ f \subset  g $.
If $x\in X$, where $ g (x)=y$, then by (\ref{EQ0097}) we have $\X \models \p [\bar x ,x]$.
Since $\qr (\p)=\qr (\f )-1\leq r-1$,
by the induction hypothesis for $\p$, $r-1$, $ g $ and $\bar x ^\smallfrown x$,
from $\X \models \p [\bar x ,x]$ it follows that $\Y \models \p [ g \bar x , g (x)]$, that is $\Y \models \p [ f \bar x ,y]$.
So for each $y\in Y$ we have $\Y \models \p [ f \bar x ,y]$, which means that $\Y \models \f [ f \bar x]$. So (\ref{EQ0093}) is true.

Now, if $\f \in \Sent _\CP$ and $\qr (\f )=r$, then, since $\Pi _r\neq \emptyset$, there is $ f \in \Pi_r$ and by (\ref{EQ0093})
(for $n=0$)  $\X \models \f $ implies $\Y \models \f$. Thus, $\X \lll _\CP \Y$.

(c) $\Rightarrow$ (a). Let $\X \lll _\CP \Y$ and, for $r<\o$, let $\Pi_r$ be the set of $ f \in \PC (\X ,\Y )$ such that there are $n\in \o$
and an $n$-tuple $\bar x \in X^n$ 
such that:

(i$_r$) $\dom  f = \{ x_0, \dots ,x_{n-1}\}$,

(ii$_r$) $\forall \f \in \CP(v_0, \dots ,v_{n-1})\;\; ( \qr (\f )\leq r \land \X \models \f [\bar x]\;\Rightarrow\; \Y \models \f [ f \bar x])$.

\noindent
First we show that $\emptyset \in I_r$, for all $r<\o$. (i$_r$) holds trivially and, since $\X \lll _\CP \Y$ we have
$\X \models \f \Rightarrow \Y \models \f $, for each $\f \in  \Sent_\CP$ (satisfying $\qr (\f )\leq r$); so (ii$_r$) holds as well.

In order to prove that the sequence $\la \Pi _r :r\in \o\ra$ satisfies (f1) and (f2), we suppose that
$ f \in \Pi_{r+1}$, $n\in \o$ and $\bar x \in X^n$, where $\dom ( f )=\{ x_0, \dots ,x_{n-1}\}$.

Since $|L|<\o$, for each $n,r<\o$ there are, up to logical equivalence, finitely many $L$-formulas $\f$
such that $|\Fv (\f )|\leq n$ and $\qr (\f )\leq r$; see \cite{}.
So, there are formulas $\p _0, \dots ,\p_{m-1}\in \CP (v_0 , \dots ,v_{n-1},v_n)$
such that $\qr (\p _j)\leq r$ and
\begin{equation}\label{EQ0094}
\forall \f \in \CP (v_0 , \dots ,v_{n-1},v_n)\;\;\Big( \qr (\f )\leq r \Rightarrow \exists j<m \;\; (\f \leftrightarrow \p _j)\Big).
\end{equation}

(f1) Let  $x\in X$ and $J:=\{ j<m : \X\models \p _j [\bar x ,x] \}$.
Then $\X\models \bigwedge _{j\in J}\p _j [\bar x ,x ]$ and, hence, $\X\models (\exists v_n \bigwedge _{j\in J}\p _j )[\bar x ]$.
Since $\qr(\exists v_n \bigwedge _{j\in J}\p _j )\leq r+1$, $\exists v_n \bigwedge _{j\in J}\p _j\in \CP$ and $ f \in \Pi_{r+1}$,
by (ii$_{r+1}$) we have $\Y\models (\exists v_n \bigwedge _{j\in J }\p _j )[ f  \bar x ]$,
so there is $y\in Y$ such that $\Y\models \bigwedge _{j\in J}\p _j [ f  \bar x ,y ]$. Thus
\begin{equation}\label{EQ0098}
\forall j\in J \;\; \Y\models \p _j [ f  \bar x ,y ]
\end{equation}
and we prove
\begin{equation}\label{EQ0095}
\forall \f \in \CP (v_0 , \dots ,v_{n-1},v_n)\;\;
\Big( \qr (\f )\leq r \land \X\models \f [\bar x ,x]\Rightarrow \Y\models \f [ f  \bar x ,y ]\Big).
\end{equation}
If $\f (\bar v,v_n)\in \CP $, $\qr (\f )\leq r$ and $\X\models \f [\bar x ,x]$
then by (\ref{EQ0094}) there is $j<m$ such that $\f \leftrightarrow \p _j$
and, hence, $\X\models \p_j [\bar x ,x]$.
So we have $j\in J$ and, by (\ref{EQ0098}), $\Y\models \p _j [ f  \bar x ,y ]$, that is $\Y\models \f [ f  \bar x ,y ]$ and (\ref{EQ0095}) is proved.

Let $ g = f \cup \{ \la x,y \ra\}$. By (\ref{EQ0095}), for each $\f (\bar v,v_n)\in \CP _0$ we have
$\X\models\f [\bar x ,x]\Rightarrow \Y\models \f [ f  \bar x ,y ]$
and, by Fact \ref{T8198'}, $ g \in \PC (\X ,\Y )$. Clearly $ f \subset  g $ and $x\in \dom ( g )=\{ x_0,\dots ,x_{n-1},x\}$ and by (\ref{EQ0095}), $ g $ satisfies (ii$_r$)
and, hence, $ g \in \Pi_r$.

(f2) Let  $y\in Y$ and let us define $J:=\{ j<m : \Y\models \neg \p _j [ f \bar x ,y] \}$ and
$\eta (\bar v):=\forall v_n \bigvee _{j\in J}\p _j(\bar v,v_n)$.
Assuming that $\X\models \eta [\bar x ]$, since $\eta \in \CP$, $\qr (\eta )\leq r+1$ and $ f \in \Pi_{r+1}$,
by (ii$_{r+1}$) we would have $\Y \models \eta [ f \bar x ]$
and, hence, for each element of $Y$ and, in particular, for $y$,
there would be $j\in J$ such that $\Y\models \p _j [ f \bar x ,y]$ which is false by the definition of $J$.
Thus $\X\models \neg \eta [\bar x ]$ and, hence, there is $x\in X$ such that
\begin{equation}\label{EQ0099}
\forall j\in J \;\; \X\models \neg \p _j [\bar x ,x ].
\end{equation}
Again  we prove (\ref{EQ0095}).
If $\f (\bar v,v_n)\in \CP $, $\qr (\f )\leq r$ and $\X\models \f [\bar x ,x]$
then by (\ref{EQ0094}) there is $j<m$ such that $\f \leftrightarrow \p _j$
and, hence, $\X\models \p_j [\bar x ,x]$.
By (\ref{EQ0099}) we have $j\not\in J$, which means that
$\Y\models \p _j [ f  \bar x ,y ]$, that is $\Y\models \f [ f  \bar x ,y ]$ and (\ref{EQ0095}) is true.

Let $ g = f \cup \{ \la x,y \ra\}$. By (\ref{EQ0095}), for each $\f (\bar v,v_n)\in \CP _0$ we have $\X\models\f [\bar x ,x]\Rightarrow \Y\models \f [ f  \bar x ,y ]$
and, by Fact \ref{T8198'}, $ g \in \PC (\X ,\Y )$. Clearly $ f \subset  g $, $y\in \ran ( g )$ and by (\ref{EQ0095}), $ g $ satisfies (ii$_r$)
and, hence, $ g \in \Pi_r$.

(a) $\Rightarrow$ (b).
Let $\X \preccurlyeq _c^{\mathrm{fin}}\Y $ and let $\la \Pi _r :r<\o\ra$ be a witness for that.
Defining $\Pi _r^{\upharpoonright } :=\{  f \in \PC (\X ,\Y ): \exists  g \in \Pi _r \;\;  f \subset  g \} $, for $r<\o$,
we show that the sequence $\la \Pi _r^{\upharpoonright }  :r<\o\ra$ satisfies (f1) and (f2).
If $ f \in \Pi _{r +1}^{\upharpoonright }$ and $ g \in \Pi _{r+1}$ where $ f \subset  g $, then, by (f1), for $x\in X$ there is $ g '\in \Pi _{r}$
such that $ g \subset  g '$ and $x\in \dom  g '$. So $ f \subset  f ':=  f \cup \{\la x,  g '(x)\ra\}\subset  g '\in \Pi _{r}$, which gives
$ f '\in \Pi _r^{\upharpoonright }$.
So (f1) is true and the proof of (f2) is similar.
Let $<_X$ and $<_Y$ be well orderings of the sets $X$ and $Y$.

By induction we prove that for each $n,r<\o$ there is a strategy $\Sigma _{n,r}$ for player II in the game $G_n^{\preccurlyeq _c}(\X ,\Y )$
such that for each play $\pi _n =\la\la x_k ,y_k \ra :k<n \ra$ we have $f_n :=\{\la x_k ,y_k \ra :k<n \}\in \Pi _r^{\upharpoonright }$.
For $n=0$ this is true, since $\emptyset \in \Pi _r^{\upharpoonright }$, for all $r<\o$.

Assuming that the statement is true for $n$ and for all $r<\o$ we define $\Sigma _{n+1,r}$ as follows.
In the first $n$ moves player II follows $\Sigma _{n,r+1}$ and obtains a partial play
$\pi _n =\la\la x_k ,y_k \ra :k<n \ra$ satisfying $f_n :=\{\la x_k ,y_k \ra :k<n \}\in \Pi _{r+1}^{\upharpoonright }$. Then

- If I chooses $x_n\in X$, then by (f1) there is $g\in \Pi _{r}^{\upharpoonright }$ such that $f_n\subset g$ and $x_n \in \dom g$;
thus $g(x_n)\in \{ y\in Y : f_n \cup \{ \la x_n ,y\ra\}\in \Pi _r^{\upharpoonright }\}$ and $\Sigma _{n+1,r}$ suggests
$
y_n =\Sigma _{n+1,r}(\pi _n ,x_n) =\min _{<_Y}\{y\in Y : f_n \cup \{ \la x_n ,y\ra\}\in \Pi _r^{\upharpoonright } \};
$

- If I chooses $y_n\in Y$, then by (f2) there is $g\in \Pi _{r}^{\upharpoonright }$ such that $f_n\subset g$ and $y_n \in \ran g$;
thus $g^{-1}(y_n)\in \{ x\in X : f_n \cup \{ \la x ,y_n\ra\}\in \Pi _r^{\upharpoonright }\}$ and $\Sigma _{n+1,r}$ suggests
$
x_n=\Sigma _{n+1,r}(\pi _n ,y_n) =\min _{<_X}\{x\in X : f_n \cup \{ \la x ,y_n\ra\}\in \Pi _r^{\upharpoonright } \}.
$

In both cases we have $f_{n+1} :=\{\la x_k ,y_k \ra :k<n+1 \}\in \Pi _{r}^{\upharpoonright }$ so our claim is true and
(b) is proved.

(b) $\Rightarrow$ (a). Suppose that for each $n\in \o$ player II has a winning strategy $\Sigma _n$ for $G_n^{\preccurlyeq _c}(\X ,\Y )$.
For $r\in \o$, let $\Pi _r$ be the set of all $f$ such that there are $j<\o$, $\bar x\in X^j$ and $\bar y\in Y^j$  such that:

(i) $f=\{ \la x_k,y_k\ra :k<j\}\in \PC (\X ,\Y)$,

(ii) There exist $m\geq r$ and a play $\pi _{j+m}=\la \la x_k,y_k\ra :k<j+m\ra$ of the game
$G_{j+m}^{\preccurlyeq _c}(\X ,\Y )$ in which II follows $\Sigma _{j+m}$.

(f1) Let $f=\{ \la x_k,y_k\ra :k<j\}\in \Pi _{r+1}$ and $x\in X$. Let $m\geq r+1$ and  let $\pi _{j+m}=\la \la x_k,y_k\ra :k<j+m\ra$ be a play of the game
$G_{j+m}^{\preccurlyeq _c}(\X ,\Y )$ in which II follows $\Sigma _{j+m}$. Then $m>0$ and in the play $\pi _{j+m}'$ in which player I
in the first $j$ moves plays in the same way as in the play $\pi _{j+m}$ and in the moves $j+1, \dots ,j+m$ plays the given $x$,
player II following $\Sigma _{j+m}$ in the moves $j+1, \dots ,j+m$ plays some $y\in Y$. So
$\pi _{j+m}'=\la \la x_0,y_0\ra , \dots ,\la x_{j-1},y_{j-1}\ra, \la x,y\ra , \dots ,\la x,y\ra \ra$, $g=f\cup \{\la x,y\ra\}\in \Pi _r$
and $x\in \dom g$ thus (f1) is true.  The proof of (f2) is similar.

(c) $\Leftrightarrow$ (d). The proof is similar to the proof of (d) $\Leftrightarrow$ (e) in Theorem \ref{T8194}.
\kdok

\section{The hierarchy of similarities}\label{S7}
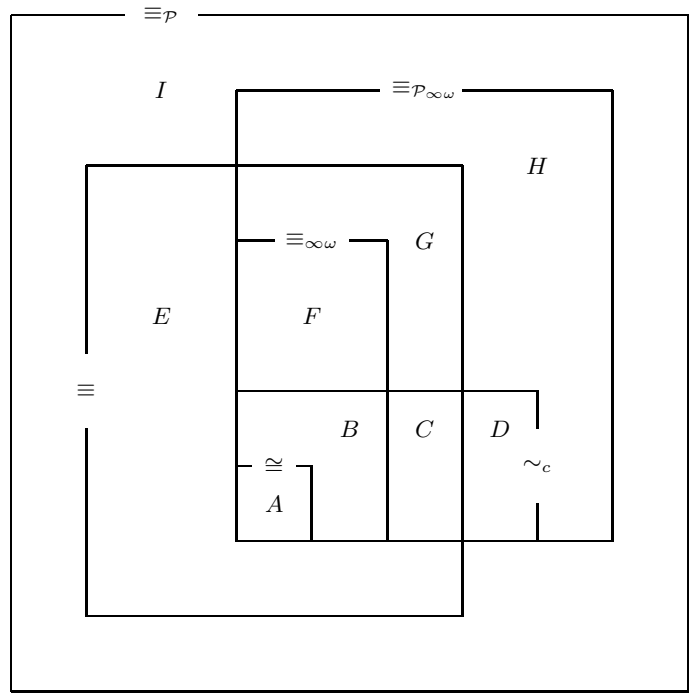
\begin{figure}[htb]
\begin{center}
\unitlength 1mm 
\linethickness{0.5pt}
\ifx\plotpoint\undefined\newsavebox{\plotpoint}\fi 

\begin{picture}(100,100)(0,0)


\put(5,5){\line(1,0){90}}
\put(15,15){\line(1,0){50}}
\put(35,25){\line(1,0){50}}
\put(35,35){\line(1,0){2}}
\put(43,35){\line(1,0){2}}
\put(35,45){\line(1,0){40}}
\put(35,65){\line(1,0){5}}
\put(50,65){\line(1,0){5}}
\put(15,75){\line(1,0){50}}
\put(35,85){\line(1,0){19}}
\put(65,85){\line(1,0){20}}
\put(5,95){\line(1,0){15}}
\put(30,95){\line(1,0){65}}
\put(5,5){\line(0,1){90}}
\put(15,15){\line(0,1){25}}
\put(15,50){\line(0,1){25}}
\put(35,25){\line(0,1){60}}
\put(45,25){\line(0,1){10}}
\put(55,25){\line(0,1){40}}
\put(65,15){\line(0,1){60}}
\put(85,25){\line(0,1){60}}
\put(95,5){\line(0,1){90}}
\put(75,40){\line(0,1){5}}
\put(75,25){\line(0,1){5}}



\footnotesize
\put(40,30){\makebox(0,0)[cc]{$A $}}
\put(40,35){\makebox(0,0)[cc]{$\cong $}}
\put(75,35){\makebox(0,0)[cc]{$\sim _c $}}
\put(60,40){\makebox(0,0)[cc]{$C $}}
\put(70,40){\makebox(0,0)[cc]{$D $}}
\put(15,45){\makebox(0,0)[cc]{$\equiv $}}
\put(25,55){\makebox(0,0)[cc]{$E $}}
\put(45,55){\makebox(0,0)[cc]{$F $}}
\put(45,65){\makebox(0,0)[cc]{$\equiv _{\infty \o} $}}
\put(60,65){\makebox(0,0)[cc]{$G $}}
\put(75,75){\makebox(0,0)[cc]{$H $}}
\put(25,85){\makebox(0,0)[cc]{$I $}}
\put(60,85){\makebox(0,0)[cc]{$\equiv _{\CP _{\infty \o}} $}}
\put(25,95){\makebox(0,0)[cc]{$\equiv _{\CP} $}}
\put(50,40){\makebox(0,0)[cc]{$B $}}
\end{picture}
\end{center}

\vspace{-5mm}
\caption{Pairs of structures \label{F0000}}
\end{figure}

In this section we compare the similarities of structures considered in this paper with the similarities $\cong$, $\equiv_{\infty \o}$ and $\equiv$ and show that the situation is as Figure \ref{F0000} describes.
We note that we construct pairs of structures of the same size and, moreover,
a pair of {\it countable} structures with the given properties, whenever it is possible.
Namely, the pairs witnessing B and F must be uncountable, because for countable structures $\X$ and $\Y$ the conditions
$\X \equiv _{\infty \o}\Y$, $\X \equiv _{\o _1 \o}\Y$ and $\X \cong \Y$ are equivalent, by Scott's theorem.
The same holds for the pairs witnessing G and H, since for countable $\X$ and $\Y$ we have $\X\equiv _{\CP _{\infty \o}} \Y$
iff $\X \sim _c\Y$; see Theorem \ref{T8194}.

Let $\CC$ be the class of structures $\X =\la X ,\r\ra$, where $|X|=\o$, $\r$ is an equivalence relation
on the set $X$ and, if $X/\r$ is the corresponding partition of $X$,

($\CC$1) $ X/\r $ contains infinitely many singletons,

($\CC$2)  For each $n\in \o$ there is an equivalence class $P\in  X/\r $ such that $|P|\geq n$.

\noindent
Defining $\CC_{\mathrm{fin}}=\{ \X \in \CC :  X/\r  \subset [X]^{<\o}\}$ and $\CC _\o=\CC \setminus \CC_{\mathrm{fin}}$
we obtain a partition $\{\CC_{\mathrm{fin}},\CC _\o  \}$ of the class $\CC$.
\begin{cla}\label{T8100}
Let $\X ,\Y \in \CC$. Then

(a) $\X \equiv _\CP \Y$;

(b) $\X \preccurlyeq _c \Y \Leftrightarrow  \X \in \CC_{\mathrm{fin}} \lor  \Y \in  \CC _\o$;

(c) $\X \sim _c \Y \Leftrightarrow \X, \Y \in \CC_{\mathrm{fin}} \lor \X, \Y \in \CC _\o$ (that is, $\CC / \sim _c =\{\CC_{\mathrm{fin}},\CC _\o  \}$).
\end{cla}
\dok
(a) Let $\X =\la X, \r\ra$ and $\Y =\la Y ,\s \ra $. For $n\in \o$
we construct a strategy $\Sigma$ for player II in the game $G_n^{\preccurlyeq _c}(\X ,\Y )$.
By ($\CC$1) we have $X':=\bigcup ((X/\r )  \cap [X]^1)\in [\o ]^{\o}$, let $Y'\in  Y/\s $, where $|Y'|\geq n$, and let
$\la X ' ,<_{X'}\ra$ and $\la Y ' ,<_{Y'}\ra$ be well orders.

Let $l<n$, let $\pi _l=\la \la x_k ,y_k \ra :k<l\ra $ be the sequence of the first $l$ moves
and $f _l=\{ \la x_k ,y_k \ra : k<l\}$.
If in the $l+1$-st move player I picks

1. $x_l=x_k$, for some $k<l$, then $y_l=\Sigma (\pi _l, x_l)=y_k$;

2. $y_l=y_k$, for some $k<l$, then $x_l=\Sigma (\pi _l, y_l)=x_k$;

3. $x_l\in X \setminus \dom f_l$, then $y_l=\Sigma (\pi _l, x_l)=\min _{<_{Y'}}(Y' \setminus \ran f_l )$;

4. $y_l\in Y \setminus \;\ran f_l\;$, then $x_l=\Sigma (\pi _l, y_l)=\min _{<_{X'}}(X' \setminus \dom f_l )$.

\noindent
Let $\pi _n =\la\la x_k ,y_k \ra :k<n\ra$ be a play of the game in which player II follows $\Sigma$.
Since in cases 3 and 4 we have $x_l \not \in \dom \f_l$ and $y_l \not \in \ran \f_l$,
case 1 ensures that $f _n=\{ \la x_k ,y_k \ra : k<n\}$ is a function and case 2 ensures that $f _n$ is an injection.
Let $\la x_k ,x_{k'}\ra\in \r $, where $k,k'<n$. If $x_k =x_{k'}$, then  $y_k =y_{k'}$ and, hence, $\la f _n (x_k ) ,f_n (x_{k'})\ra\in \s $.
If $x_k \neq x_{k'}$, then, since $\la x_k ,x_{k'}\ra\in \r $, we have $x_k , x_{k'}\not\in X'$, which means that,
at the step when the pairs $\la x_k ,y_k \ra$ and $\la x_{k'} ,y_{k'} \ra$ were chosen for the first time we had case 3;
so $y_k ,y_{k'}\in Y' $ and, hence, $\la f _n (x_k ) ,f_n (x_{k'})\ra =\la y_k ,y_{k'}\ra \in \s $.
Thus $f _n\in \PC (\X ,\Y )$ and player II wins.

(b) ($\Rightarrow$) On the contrary, suppose that $\X \preccurlyeq _c \Y$, $X' \in  ( X/\r ) \cap [X]^{\o}$ and $ Y/\s  \subset [Y]^{<\o}$.
Then there would be $f\in \Cond (\X ,\Y)$ and, hence, $f[X']\subset Y'$, for some $Y'\in  Y/\s $, which is impossible.

($\Leftarrow$)
Let $ X/\r  =\{ \{s_n \} :n< \o\}\cup \{ X_i: i\in \o\}$, where $|X_i|>1$, for $i<\o$.
Then by ($\CC$2) we have $|\bigcup _{i<\o}X_i|=\o$.

First, if there exists $Y'\in  (Y/\s )  \cap [Y]^\o$ and $f:X\rightarrow Y$ is a bijection satisfying $f[\bigcup _{i<\o}X_i]=Y'$
and $f[\{s_n: n<\o\}]=Y\setminus Y'$, then $f\in \Cond (\X ,\Y)$.

Second, suppose that $|X_i|<\o$, for all $i<\o$ and that $ Y/\s  \subset [Y]^{<\o}$.
Let $ Y/\s  =\{ \{t_n \} :n< \o\}\cup \{ Y_j: j\in \o\}$, where $1< |Y_j|<\o$, for $j<\o$.
By recursion we define a sequence $\la j_i : i<\o\ra$ such that for each $i<\o$ we have:

(i) $|X_i|\leq |Y_{j_i}|$,

(ii) $i'< i \Rightarrow j_{i'} \neq j_i$.

\noindent
First, by ($\CC$2) there is $j_0<\o$ such that $|X_0|\leq |Y_{j_0}|$. Let $i<\o$ and suppose that
$\la j_{i'} : i'<i \ra$  is a sequence satisfying (i) and (ii). By ($\CC$2) we can choose (minimal) $j_i<\o$ such that
$|Y_{j_i}|>\max \{|X_i|, |Y_{j_0}|, \dots , |Y_{j_{i-1}}|\}$ and the sequence $\la j_{i'} : i'\leq i \ra$ satisfies (i) and (ii).
The recursion works.

By (i) there are injections $f_i : X_i \rightarrow Y_{j_i}$, for $i<\o$, and, since by (ii) the sets $Y_{j_i}$, $i<\o$,
are pairwise disjoint, $f:= \bigcup _{i<\o }f_i $ is an injection which maps the set $X\setminus \{ s_n :n\in \o\}$
onto a subset $Y'$ of  $\bigcup _{i<\o }Y_{j_i}$, which means that $|Y\setminus Y'|=\o$.
So there is a bijection $g:\{ s_n :n\in \o\}\rightarrow Y\setminus Y'$ and, clearly, $f\cup g\in \Cond (\X ,\Y)$.

(c) By (b), if $\X, \Y \in \CC_{\mathrm{fin}}$ or $\X, \Y \in \CC _\o$, then we have $\X \sim _c \Y$.
On the other hand, if $\X \in \CC _\o$ and $\Y \in \CC_{\mathrm{fin}}$, then, by (b) again, $\X\not\preccurlyeq _c \Y$ and, hence, $\X\not\sim _c \Y$
\kdok
\begin{fac}\label{T8199}
If $\X$ is a countable equivalence relation, there is a first-order theory $\CT _\X^{\mathrm{fin}}$
determining the number ($\leq \o$) of $k$-sized equivalence classes, for each $k\in \N$.

If, in addition,  for each $n\in \N$ there is a finite equivalence class of size $\geq n$, then the theory $\CT _\X^{\mathrm{fin}}$
is complete but not $\o$-categorical.
\end{fac}
The following examples show that the situation is as Figure \ref{F0000} describes.
\begin{ex}\label{EX0001}\rm
A: Of course here we can take $\X \cong \X$, for any structure $\X$.

B: $\X \sim _c \Y$, $\X \equiv _{\infty \o}\Y$, and $\X \not\cong \Y$.
Let $\{ A, B\}\cup \{ C_\a :\a <\o _1\}$ be a partition of $\o _1$, where
$A,B \in [\o _1]^{\o _1}$ and $C_\a\in [\o _1]^\o$, for $\a <\o _1$,
and let $\X =\la \o _1 ,\r \ra$ and $\Y =\la \o _1 ,\s \ra$, where $\r$ and $\s$
are the equivalence relations on $\o _1$ corresponding to the partitions $\{ A\cup B\}\cup \{ C_\a :\a <\o _1\}$ and
$\{ A, B\}\cup \{ C_\a :\a <\o _1\}$, respectively.
It is easy to see that $\X \sim _c \Y$ and $\X \not\cong \Y$ (see Example 3.13 from \cite{KuMo1}).
In addition if $V[G]$ is a generic extension of the universe by the collapsing algebra $\Col (\o _1 ,\o)$, then, in $V[G]$,
$\X$ and $\Y$ are countable structures with one equivalence relation, having $\o$-many equivalence classes and all of them are infinite.
Thus $V[G]\models \X \cong \Y$ and, hence,  $\X \equiv _{\infty \o}\Y$.

E: $\X \equiv \Y$ and $\X \not\equiv _{\CP _{\infty \o}}\Y$. Let $\X$ and $\Y$ be linear orders, where
$\X \cong \o$ and $\Y \cong\o + \zeta$ ($\zeta =\o ^* +\o$ is the order type of the integers).
It is well known that $\X \equiv \Y$. Assuming that $\X \equiv _{\CP _{\infty \o}}\Y$, by Theorem \ref{T8194}
we would have $\X \sim _c \Y$ and, since all linear orders are reversible and, hence, weakly reversible, $\X \cong \Y$,
which is false.

D: $\X \sim _c \Y$ and $\X \not\equiv \Y$. Let $\X =\la X, \r\ra$ be the Rado graph (see \cite{Camer}), let $x,y\in X$, where
$x\,\r \,y$, that is, $\la x, y\ra , \la y,x \ra \in \r$ and let $\Y =\la X , \s \ra$, where $\s =\r \setminus \{ \la x,y\ra\}$.
Then, since the relation $\s$ is not symmetric, we have $\X \not\equiv \Y$.
In addition, the structure $\X _1=\la X , \r _1 \ra$, where $\r _1 =\s \setminus \{ \la y,x\ra\}$
is a graph obtained from $\X$ by deleting an edge and $\X _1\cong \X$ (see \cite{Camer}). Now we have
$\r \cong \r _1 \subset \s \subset \r$, which implies that $\X \preccurlyeq _c \Y \preccurlyeq _c \X$ and, thus, $\X \sim _c \Y$.

I: $\X \equiv_\CP \Y$, $\X \not\equiv_{\CP _{\infty \o}} \Y$  and $\X \not\equiv \Y$. Let $\X \in \CC_{\mathrm{fin}}$ and $ \Y \in \CC _\o$
(see the beginning of this section), where $ X/\r  =\{ \{s_n \} :n< \o\}\cup \{ X_i: i\geq 2\}$ and $|X_i|=i$, for $i\geq 2$,
and $ Y/\s  =\{ \{t_n \} :n< \o\}\cup \{ Y'\}$, where $|Y'|=\o$. By Claim \ref{T8100}(a) we have  $\X \equiv_\CP \Y$
and by (c) we have  $\X \not\sim_c \Y$, which, together with Theorem \ref{T8194} gives $\X \not\equiv_{\CP _{\infty \o}}\Y$.
For example we have $\X \models \neg \f$ and $\Y \models  \f$ where the sentence $\f \in \CP _{\infty \o}$ says that there is an infinite equivalence class:
$$\textstyle
\f :=\exists v \bigwedge _{n\in \N} \exists v_0 , \dots ,v_{n-1} \Big(\bigwedge _{i<j<n} v_i\neq v_j \land \bigwedge _{i<n}R(v,v_i)\Big).
$$
Finally, we have $\X \not\equiv \Y$, because $\X \models \p$ and $\Y \models \neg \p$ where $\p$ is the first order $L_b$-sentence
saying that there is an equivalence class of size 2,
$$\textstyle
\p :=\exists u,v \;\;\Big( u\neq v \land R(u,v) \land \forall w \;\;(w\neq u \land w\neq v \Rightarrow \neg R(w,u))\Big).
$$

F: $\X \equiv_{\infty \o} \Y$ and $\X \not\sim _c \Y$. Let $\X$ and $\Y$ be the real and the irrational line respectively.
Since linear orders are reversible, $\X \sim _c \Y$ would imply that $\X\cong \Y$, which is false; so, $\X \not\sim _c \Y$.
Further, in a generic extension  $V_{\Col (2^\o ,\o)}[G]$,
$\X$ and $\Y$ are countable dense linear orders without end points and, by Cantor's theorem,
$V[G]\models \X \cong \Y \cong \Q$ or, equivalently,  $\X \equiv _{\infty \o}\Y$.

H: $\X \equiv_{\CP _{\infty \o}} \Y$, $\X \not\equiv \Y$ and $\X \not\sim _c \Y$.
Let $\X =\la \o _1 ,\r \ra$ and $\Y =\la \o _1 ,\s \ra$, where $\r$ and $\s$
are the equivalence relations on $\o _1$ corresponding to the partitions
$\{ 2, \o \setminus 2 \}\cup \{ \{\a \} :\o \leq \a <\o _1\}$ and
$\{ \{ n\} : n\in \o \}\cup \{ \o _1 \setminus \o\}$, respectively.
Then $\X \not\equiv \Y$, because $\X$ has an equivalence class of size 2, which is false in $\Y$.
Since $\Y$ has an equivalence class of size $\o _1$, while all equivalence classes of $\X$ are at most countable,
we have $\Y \not\preccurlyeq _c \X$ and, hence, $\X \not\sim _c \Y$.
Finally, in a generic extension  $V_{\Col (\o_1 ,\o)}[G]$,
$\X$ and $\Y$ are countable structures belonging to the class $\CC _\o$
and by Claim \ref{T8100}(c) we have $\X \sim _c \Y$.
Thus $V[G]\models \X \sim _c \Y $ and, by Theorem \ref{T8194}, $\X \equiv_{\CP _{\infty \o}} \Y$.

C: $\X \equiv \Y$,  $\X \sim _c \Y$ and $\X \not\equiv_{\infty \o} \Y$.
Let $\f$, $\p _n$, for $n\geq 2$, and $\eta _n$,  for $n\in \N$, be $L_b$-sentences such that
$\f$ says that $R(u,v)$ is an equivalence relation,
$\p _n$ says that there is exactly one equivalence class of size $n$,
and $\eta _n$ says that there are at least $n$ equivalence classes of size 1.
Then, by Fact \ref{T8199}, the theory $\CT =\{\f \} \cup \{\p _n :n\geq 2\} \cup \{\eta _n :n\in \N\}$ is complete.
Let $\X $ and $\Y$ be countable models of $\CT$, where $\X$ has one infinite equivalence class,
while $\Y$ has two. Then $\X \equiv \Y$ and, since $\X ,\Y \in \CC _\o$,  by Claim \ref{T8100}(c) we have $\X \sim _c \Y$.
But $\X \not\cong \Y$ which by Scott's theorem implies $\X \not\equiv_{\o _1\o} \Y$ and, hence, $\X \not\equiv_{\infty \o} \Y$.

G: $\X \equiv \Y$, $\X \equiv_{\CP _{\infty \o}} \Y$, $\X \not\sim _c \Y$ and $\X \not\equiv_{\infty \o} \Y$. Let $\CT$ be the theory from part C
and let $\X $ and $\Y$ be models of $\CT$ of size $\o _1$ having countably many equivalence classes of size 1 and
such that $ (X/\r ) \setminus [X]^{<\o}=\{ X'\}$ and $ (Y/\s ) \setminus [Y]^{<\o}=\{ Y',Y''\}$, where
$|X'|=|Y'|=|Y''|=\o _1$. Then $\X \equiv \Y$. If $f$ is a homomorphism from $\X$ to $\Y$, then $f$ must map
$X'$ into an equivalence class of $\Y$ and, hence, $f[X']\subset Y'$ or $f[X']\subset Y''$; thus $f$ can not be
a surjection, which shows that $\X \not \preccurlyeq _c \Y$ and, consequently, $\X \not\sim _c \Y$.
Again, in a generic extension  $V_{\Col (\o_1 ,\o)}[G]$,
$\X$ and $\Y$ are countable structures belonging to the class $\CC _\o$
and by Claim \ref{T8100}(c) we have $\X \sim _c \Y$.
Thus $V[G]\models \X \sim _c \Y $ and, by Theorem \ref{T8194}, $\X \equiv_{\CP _{\infty \o}} \Y$.
But $V[G]\models \X \not\cong \Y $ and, thus, $\X \not\equiv_{\infty \o} \Y$.
\end{ex}
\noindent
{\bf Acknowledgments.}
This research was supported by the Ministry of Education and Science of the Republic of Serbia (Project 174006).

\end{document}